# Some consequences of a recursive number-theoretic relation that is not the standard interpretation of any of its formal representations

Bhupinder Singh Anand

We give precise definitions of primitive and formal mathematical objects, and show: there is an elementary, recursive, number-theoretic relation that is not a formal mathematical object in Gödel's formal system P, since it is not the standard interpretation of any of its representations in P; the range of a recursive number-theoretic function does not always define a formal mathematical object (recursively enumerable set) consistently in any Axiomatic Set Theory that is a model for P; there is no P-formula, *Con*(P), whose standard interpretation is unambiguously equivalent to Gödel's number-theoretic definition of "P is consistent"; every recursive number-theoretic function is not strongly representable in P; Tarski's definitions of "satisfiability" and "truth" can be made constructive, and intuitionistically unobjectionable, by reformulating Church's Thesis constructively; the classical definition of Turing machines can be extended to include self-terminating, converging, and oscillating routines; a constructive Church's Thesis implies, firstly, that every partial recursive number-theoretic function has a unique, constructive, extension as a total function, and, secondly, that we can define effectively computable number-theoretic functions that are not classically Turing-computable; Turing's and Cantor's diagonal arguments do not necessarily define Cauchy sequences.

# Contents[1]

## I. A meta-theorem of recursive asymmetry

## II. Consequences

---

[1] *Key words*: Algorithm, arithmetic, Cantor, Cauchy, Church, computable, consistent, constructive, Dedekind, diagonal construction, effective method, expressible, first order, formal system, function, general recursive, Godel, Hilbert, individually terminating, interpretation, intuitive, mapping, mathematical object, model, natural number, neo-classical, non-recursive, number-theoretic, oscillating, partial recursive, Peano, predicate calculus, primitive recursive, proposition, P vs NP, real number, recursive, recursively enumerable, replacement axiom, representable, satisfiable, self terminating, sentence, set theory, sound, standard model, Tarski, total function, truth, Turing, uncomputable, uncountable, undecidable, uniformly terminating

# I.  A meta-theorem of recursive asymmetry

## 1.  Introduction

### 1.1  Preamble

In this paper, we essentially address the question: Are the concepts "non-algorithmic"[2]

and "non-constructive"[3] necessarily synonymous in classical[4] logic and mathematics?

We consider, as a natural starting point, the classical interpretation of the reasoning and

conclusions in Gödel's seminal 1931 paper [Go31a]. Gödel argues that, in a formal

language[5] as basic as Peano Arithmetic[6], there are undecidable[7] sentences[8] that can be

---

[2] We follow Mendelson's definition of an algorithm as an effectively computable function ([Me64], p208). Subject to their being individually proven as formal mathematical objects - a concept that we define precisely - we also follow Mendelson's set-theoretic definitions of a "function" and of a "relation" ([Me64], p6-7). In other words, since we argue later that we can define number-theoretic functions that are not definable as formal mathematical objects in any Axiomatic Set Theory, we treat the sets in Mendelson's definitions as hypothetical and intuitive, but not formal, mathematical objects. Therefore, assuming formal set-theoretic properties for them, even in informal reasoning, may invite inconsistency.

[3] We note that the term "constructive" - and its synonym "effective" - is used both in its familiar linguistic sense, and in a mathematically precise sense. Mathematically, we term a concept as "constructive" if, and only if, it can be defined in terms of pre-existing concepts without inviting inconsistency (cf. Mendelson's remarks in (cf. Mendelson's remarks in [Me64], p82). Otherwise, we understand it in an intuitive sense to mean unambiguously verifiable, by some "effective method" ([Me64], p207-8), within some finite, well-defined, language or meta-language ([Me64], p31, footnote). Generally, it may be taken to correspond, broadly, to Gödel's concept of "intuitionistically unobjectionable" ([Go31a], p26).

[4] For the purposes of this paper, we take the expositions by Hardy [Ha47], Landau [La51], Mendelson [Me64], Rudin [Ru53] and Titchmarsh [Ti61] as standard presentations of classical mathematical reasoning and conclusions.

[5] By a "formal language" we mean a "formal system" or a "formal theory" as in Mendelson ([Me64], p29).

[6] By Peano Arithmetic, we mean any formalisation of Dedekind's formulation of his Peano Postulates (cf. [Me64], p102). We take Gödel's formal system P ([Me64], p9) as a second order formalisation of Peano

recognised as true[9] under classical interpretation[10], but which are not provable[11] within

the system. Does this imply that such recognition, in some cases, cannot be duplicated in

any artificially constructed and, more important, artificially controlled, mechanism or

organism whose design is based on classical logic?[12]

The scientific, and philosophical, dimensions of an affirmative answer to the last question have been broadly reviewed, and addressed, by Roger Penrose in [Pe90] and [Pe94].

Penrose's argument is based on a strongly Platonist thesis that sensory perceptions

---

Arithmetic, and Mendelson's formal system S ([Me64], p102) as a standard formalisation of classical first order Peano Arithmetic.

[7] We follow Mendelson's definition of an "undecidable sentence" ([Me64], p143).

[8] When referring to a formal language, we assume the terms "sentence" and "proposition" are synonymous, and that they refer to a well-formed expression of the language that contains no free variables, and which translates, under an interpretation, as a proposition in the usual, linguistic, sense.

[9] We note that the term "true" is used both in its familiar linguistic sense, and in a mathematically precise sense; the appropriate meaning is usually obvious from the context. Mathematically, we follow Mendelson's exposition of the truth of a formal sentence under an interpretation as determined by Tarski's definitions of satisfiability and truth ([Me64], p51).

[10] We note that the term "interpretation" is also used both in its familiar linguistic sense, and in a mathematical sense; the appropriate meaning is usually obvious from the context. Mathematically, we follow Mendelson's definitions of "interpretation" ([Me64], §2, p49), and of "standard interpretation" ([Me64], p107). We note that the interpreted relation $R(x)$ is obtained from the formula $[R(x)]$ of a formal system P by replacing every primitive, undefined symbol of P in the formula $[R(x)]$ by an interpreted mathematical symbol (i.e. a symbol that is a shorthand notation for some, semantically well-defined, concept of classical mathematics).

So the P-formula $[(Ax)R(x)]$ interprets as the sentence $(Ax)R(x)$, and the P-formula $[\sim(Ax)R(x)]$ as the sentence $\sim(Ax)R(x)$.

We also note that the meta-assertions "$[(Ax)R(x)]$ is a true sentence under the interpretation M of P", and "$(Ax)R(x)$ is a true sentence of the interpretation M of P", are equivalent to the meta-assertion "$R(x)$ is satisfied for any given value of $x$ in the domain of the interpretation M of P" ([Me64], p51).

[11] The term "provable", when applied to a formal expression $F$, implies the existence of a finite sequence of formal expressions (of which $F$ is the last) that is constructed by a finite set of rules, and which forms the basis of Mendelson's definition of a "proof" ([Me64], p29).

[12] We note that the question may have economic significance globally, particularly in areas relating to the development of strategic and infra-structural products, facilities and services that are based on the proposed replication of human intelligence by artificial mechanisms or organisms.

simply mirror aspects of a universe that exists, and will continue to exist, independent of any observer ([Pe90], p123, p146)[13]. On this view, individual consciousness would be a discovery of what there is (cf. [Pe90], p124), and be independent of the language in which such discovery is expressed. It follows that recognition of intuitive truth would be individually asserted - and, implicitly, fallible - correlations between the unverifiable - and, ipso facto, infallible - intuitive experiences of an individual consciousness, and the formal expressions of a communicable language.

The issue, then, is whether classical logic can adequately formalise intuitive truth to make it infallible, or whether such recognition is essentially fallible[14]. Penrose opts for the inadequacy of classical logic to completely capture a Platonic mathematical reality that, he believes, manifests itself, firstly, in thought - which originates in the mind consequent to sensory experience - and, secondly, in its representation in a communicable language. He supports his view by highlighting the "ethereal" presence, and non-verifiable properties, of various non-algorithmic ([Pe90], p168), and implicitly non-constructive, mathematical concepts that are accepted in our formal languages as essential to classical mathematics ([Pe90], p123-8).

Although Penrose's arguments represent only one, and perhaps an arguably extreme, point of view[15], they emphasise that classical mathematics may yet need to adequately legitimise the acceptance, into a theory, of formally definable mathematical objects (cf. [Pe90], p147)[16], most obviously those that can be argued as being essentially non-constructive.

---

[13] An obvious, but arguably relevant, objection to this argument is that it assumes multiple, spatially separated, observers can each, Deity-like, acquire identical knowledge of a Universe simultaneously without altering, or even indirectly influencing, the knowledge that is sought to be acquired.

[14] Of course, there is an inescapable element of circularity in considering the fallibility of assertions that are asserted as intuitively true.

[15] See Psyche, Vol. 2(9), June 1995, Symposium on Roger Penrose's Shadows of the Mind.

[16] We give a precise definition of the term "mathematical object" in a later paragraph.

Now we note that Penrose appears to base his thesis on, amongst others, a classical consequence of Gödel's reasoning and conclusions: we cannot express Tarskian definitions[17] of the "satisfiability" and "truth" of formal expressions under an interpretation algorithmically ([Pe90], p159)[18]. He concludes from this that, although we may follow a common intuitive process for discovering common mathematical aspects of the universe, not all our mathematically expressible discoveries are expressible by classical algorithms ([Pe90], p533, p548).

However, Penrose's arguments also appear to imply further, albeit implicitly, that our recognition of intuitive "arithmetical truth" - even when this is accepted as being adequately formalised by the classical Tarskian definitions of the "satisfiability" and "truth" of formal expressions under an interpretation - is "absolutely" non-constructive (cf. [Pe90], p145-6).

Thus, although Penrose does not seem to question the mathematical form of Church's Thesis[19] ([Pe90], p64-65) - which, essentially, postulates that every effectively computable function is algorithmic - he seems to conclude from his arguments, concerning the inadequacy of classical logic, that there may be non-algorithmic, non-constructive, ways of acquiring mathematical insight and knowledge ([Pe90], p538). As is evidenced in his discussion of Lucas' Gödelian argument ([Pe90], p539), Penrose does not appear to entertain the possibility that there may be non-algorithmic reasoning that could be intuitionistically accepted as constructive; his arguments seem to, implicitly, treat the terms "non-algorithmic" and "non-constructive" as synonymous[20].

---

[17] We take Mendelson's exposition ([Me64], p49-52) as representative of the classical Tarskian definitions of the "satisfiability" and "truth" of well-formed formulas of a formal language under a given interpretation.

[18] This is, essentially, an intuitive interpretation of Tarski's Theorem: The set *Tr* of Gödel-numbers of the formal expressions of a first order Peano Arithmetic that are true in the standard model is not arithmetical ([Me64], p151).

[19] We take the classical Church Thesis as the assertion that a number-theoretic function is effectively computable (partially) if, and only if, it is (partially) recursive ([Me64], p147, p227).

## 1.2 A surprising theorem

In this paper, however, we argue that what Penrose, for instance, views as the essentially, and absolutely non-constructive, aspects of mathematical concepts, may simply be manifestations of a *removable* ambiguity in the classical Tarskian definitions of the satisfiability, and truth, of our formal expressions under an interpretation.

We argue that this leads to an alternative interpretation of the classical consequences of Gödel's seminal 1931 paper [Go31a], with implications for the foundations of philosophy, logic, mathematics and computability.

We start by noting that the introduction of classical, non-constructive, objects into our mathematical ontology - particularly those introduced through unrestricted definitions[21] - can be qualified by formally defining a set as a mathematical object:

---

[20] Professor Martin Davis, Courant Institute of Mathematical Sciences, New York, critically reviews this particular aspect of Penrose's argument in an article dated 22nd March 1995, titled "Is Mathematical Insight Algorithmic?", that the author downloaded from an unrecorded source on the web. He argues that: "... Gödel's incompleteness theorem (in a strengthened form based on work of J.B. Rosser as well as the solution of Hilbert's tenth problem) may be stated as follows: There is an algorithm which, given any consistent set of axioms, will output a polynomial equation P = 0 which in fact has no integer solutions, but such that this fact can not be deduced from the given axioms. Here then is the true but unprovable Gödel sentence on which Penrose relies and in a particularly simple form at that. Note that the sentence is provided by an algorithm. If insight is involved, it must be in convincing oneself that the given axioms are indeed consistent, since otherwise we will have no reason to believe that the Gödel sentence is true". As we argue in this paper, however, the real issue is not whether there is an algorithm that outputs P = 0, but whether, for any given set of natural number values for its free variables, the fact that P = 0 has no integer solutions can be determined in a non-algorithmic, yet constructive way.

[21] These would include, for instance, definitions such as those of: a Turing-uncomputable number-theoretic function ([Tu36], 9, para II); a set-theoretical limit of an iterative process in the completion of a metric space, such as the Cantor set ([Ru53], p34); a non-recursive recursively-enumerable set (cf. [Me64], p251, Proposition 5.20(5)).

**Definition 1**(*i*): A *primitive mathematical object* is any symbol for an individual constant, predicate letter, or a function letter (cf. [Me64], p46; also p1, p10), which is defined as a primitive symbol of a formal mathematical language.

**Definition 1**(*ii*): A *formal mathematical object* is any symbol for an individual constant, predicate letter, or a function letter that is either a primitive mathematical object, or that can be introduced through definition (cf. [Me64], p82) into a formal mathematical language without inviting inconsistency[22].

**Definition 1**(*iii*): A *mathematical object* is any symbol that is either a *primitive mathematical object*, or a *formal mathematical object*.

**Definition 1**(*iv*): A *set* is the range of any function whose function letter is a mathematical object.

The significance of these definitions is seen in Meta-theorem 1. We prove, there, the existence of an asymmetrical recursive number-theoretic relation[23] - one that is intuitively decidable constructively, but which cannot be introduced through definition as a formal mathematical object into any formal system of Peano Arithmetic without inviting inconsistency; nor, ipso facto, into any Axiomatic Set Theory that models[24] (cf. [Me64], p192-3) such Arithmetic. Hence, it is not a formal mathematical object, and the range of its characteristic function[25] is not a recursively enumerable set[26]!

---

[22] We take Mendelson's Corollary 1.15 ([Me64], p37), as the classical meta-definition of consistency.

[23] We treat the terms "relation" and "predicate" as synonyms, and use them interchangeably.

[24] We follow Mendelson's definition of a model ([Me64], p51): An interpretation is said to be a model for a set $T$ of well-formed formulas of P if, and only if, every well-formed formula in $T$ is true for the interpretation.

[25] If $R(x)$ is a relation (predicate), then the characteristic function $C_R(x)$ is defined as follows: $C_R(x) = 0$ if $R(x)$ is true, and $C_R(x) = 1$ if $R(x)$ is false (5, p119).

This is an astonishing result[27]! Firstly, recursive number-theoretic functions and relations are classically accepted as the basic building blocks for defining constructive, and intuitionistically unobjectionable, mathematical objects[28]. Secondly, and in vivid contrast, the relative consistency, and independence, of the Continuum Hypothesis would imply, prima facie, that we may also treat Cantor's non-constructive cardinal, *aleph*₁, as a valid formal mathematical object[29]; thus, we may introduce axiomatic definitions - and an individual constant symbol - for it into any Axiomatic Set Theory without inviting inconsistency!

### 1.3  An alternative interpretation of Gödel's Proof

Now we note that, in the proof of Theorem VI of his 1931 paper ([Go31a], p24), Gödel argues that, in any consistent, formal, system P that formulates Peano's Arithmetic, we can construct a valid expression of the system, say $[R(x)]$[30], such that $[R(n)]$ is P-provable for any given numeral $[n]$[31], but $[(Ax)R(x)]$ is not P-provable. The classical interpretation of this is that although $[(Ax)R(x)]$ is not P-provable, it is true under its standard interpretation by Tarski's definitions of satisfiability and truth ([Me64], p51).

---

[26] A recursively enumerable set is classically defined as the range of some recursive number-theoretic function, and is implicitly assumed consistent with any Axiomatic Set Theory that is a model for P ([Me64], p250).

[27] Loosely speaking, it may be viewed as a constructive arithmetical parallel to Russell's non-constructive, and paradoxical, set ([Me64], p2). Some philosophical implications of this result are also echoed in David Chalmers remark: "... I have some sympathy with Penrose's idea that we have an underlying sound competence, even if our performance sometimes goes astray. But further, it seems to me that to hold that this is the only problem in Penrose's argument would be to concede too much power to the argument. It would follow, for example, that there are parts of our arithmetical competence that no *sound* formal system could ever duplicate; it would seem that our unsoundness would be essential to our capacity to see the truth of Gödel sentences, for example. This would be a remarkably strong conclusion, and does not seem at all plausible to me". Minds, Machines, And Mathematics, para 2.5, Psyche, Vol. 2(9), June 1995.

[28] See, for instance, Gödel's remarks ([Go31a], p23, footnote 39).

[29] This, essentially, seems to reflect Gödel's point of view, which he expresses in his 1947 paper, "What is Cantor's continuum problem?", whilst discussing whether Cantor's continuum hypothesis should be added to set theory as a new axiom. [Kurt Gödel, Collected Works, vol. 2, Oxford University Press, 1986–2003.]

[30] We use square brackets to differentiate between a formal expression $[R(x)]$ and its interpretation $R(x)$.

[31] We follow Gödel's definition of a numeral ([Go31a], p10).

We argue, however, that, by implications that are implicit in Tarski's definitions, $[R(n)]$ may be viewed alternatively as an expression whose standard interpretation, $R(n)$, can be effectively asserted as holding *individually* - and not necessarily algorithmically - for any given natural number $n$, but $R(x)$ cannot be effectively asserted as holding *uniformly* - in the sense of algorithmically - for all natural numbers $x$ collectively.

## 1.4  The ambiguity

Thus, we deny the very basis for the above interpretation of Penrose's thesis, and argue that classical interpretations of Tarski's definitions of satisfiability and truth contain an ambiguity: they implicitly imply the existence of an ambiguous *effective method*[32] for deciding whether formal expressions such as $[R(x)]$ are satisfied under a given interpretation. Specifically, they fail to entertain the possibility that such a method may be non-algorithmic[33]. In other words, for any given value $n$ of its free variable under a given interpretation, there may always be a - possibly $n$-specific - (non-algorithmic) effective method that can effectively decide whether the interpretation $R(n)$ of a formal expression such as $[R(n)]$ holds *individually*, even when there is no $n$-independent (algorithmic) effective method  that can effectively decide whether the expression $[R(x)]$ is satisfied *uniformly*, under a given interpretation, when we substitute any numeral $[n]$ for its free variable.

---

[32] The phrase "effective method" (as well its synonym, "mechanical procedure") is not at all precise; as Mendelson notes ([Me64], p207), "... what we mean is a process which requires no ingenuity for its performance".

[33] In other words, we argue that not every effective method is necessarily algorithmic, although every algorithm is an effective method. The possibility that "truth" may be non-algorithmic, and yet constructive, is implicit in one of Gödel's unpublished essays, "Some basic theorems on the foundations" (Kurt Gödel, Collected Works, vol. 3, Oxford University Press, 1986–2003.): "I wish to point out that one may conjecture the truth of a universal proposition (for example, that I shall be able to verify a certain property for any integer given to me) and at the same time conjecture that no general proof for this fact exists. It is easy to imagine situations in which both these conjectures would be very well founded. For the first half of it, this would, for example, be the case if the proposition in question were some equation F(n) = G(n) of two number-theoretical functions which could be verified up to very great numbers n". It is also implicit in Turing's remarks ([Tu36], §9, para II).

We argue, further, that the above ambiguity is removable by making the above possibility explicit in the classical Tarskian definitions of satisfiability and truth, and by introducing a constructive expression of Church's Thesis. We then argue, for instance, that every classically Turing-computable number-theoretic function is effectively computable individually. Consideration of the converse leads naturally to the definition of self-terminating, converging and oscillating neo-classical Turing routines, and an argument that the Church-Turing Thesis may be false, since we can now define effectively computable functions that are not classically Turing-computable. We also show that Turing's "uncomputable" real numbers, and Cantor's "uncountable" real numbers, do not necessarily correspond to Cauchy sequences of rational numbers; thus, they cannot be assumed to always define Dedekind real numbers.

## 1.5 Overview

In his 1931 paper [Go31a], Gödel prefaces his Theorem V with the remark ([Go31a], p22):

> "The fact which can be vaguely formulated as the assertion that every recursive relation is definable within the system P (under its intuitive interpretation), is rigorously expressed by the following theorem, without reference to the intuitive meaning of the formulas of P.
>
> Theorem V: For every recursive relation $R(x_1, ..., x_n)$, there is an $n$-ary PREDICATE $r$ (with the FREE VARIABLES $u_1 ... u_n$) such that, for all $n$-tuples of numbers $(x_1, ..., x_n)$, we have:
>
> $$R(x_1, ..., x_n) => \text{Bew}[Sb(r\,(u_1 ... u_n)|(Z(x_1) ... Z(x_n)))] \qquad (3)$$

$\sim R(x_1, ..., x_n) \Rightarrow \text{Bew}[\text{Neg } Sb(r\,(u_1 ... u_n)|(Z(x_1) ... Z(x_n)))]$      (4)"

In a footnote, he adds that ([Go31a], footnote 38):

"The VARIABLES $u_1 ... u_n$ can be arbitrarily prescribed. There always exists, e.g. some $r$ with the FREE VARIABLES 17, 19, 23, etc., for which (3) and (4) hold."

He then qualifies his proof with the remark ([Go31a], p23):

"We shall be content here to indicate the outline of the proof of this theorem, since it offers no theoretical difficulties and is fairly tedious."

In another footnote, he clarifies that ([Go31a], footnote 39):

"Theorem V depends of course upon the fact that, for a recursive relation $R$, it is decidable on the basis of the axioms of the system P whether or not $R$ holds for any given $n$-tuple of numbers."

Since Gödel does not give a rigorous proof of the theorem, it is not clear whether his remark - that the recursive relation $R(x_1, ..., x_n)$ is "... definable within the system P (under its intuitive interpretation) ..." - means that his reasoning is based on an implicit assumption that $R(x_1, ..., x_n)$ is the standard interpretation of some PREDICATE $r$ of P. Such an assumption would imply that a predicate letter for the $n$-ary relation "$R$", along with suitable defining axioms, could be introduced into P (cf. [Me64], p82) without inviting inconsistency. We address this issue, and its possible direct, and indirect, consequences, in the following sections.

In Meta-theorem 1, we consider the elementary, recursive, number-theoretic relation $x=Sb(y\ 19|Z(y))$ (cf. 1, p20, def. 31), and prove that it is not the standard interpretation of any[34] of its formal representations[35] in Gödel's system P.

In Meta-lemma 1, we then conclude that we cannot introduce a finite number of arbitrary recursive number-theoretic functions and relations, as function letters and predicate letters respectively, into a formal system of Arithmetic, such as P, without inviting inconsistency. We thus conclude, in Corollary 1.1, that every recursive number-theoretic relation does not, consistently, well-define a (recursively enumerable) sub-set of the natural numbers in any Axiomatic Set Theory[36] that is a model for P; and, in Corollary 1.2, that not every constructively well-defined number-theoretic function is a mathematical object (so it may not define a mapping[37] in any Axiomatic Set Theory that is a model for P).

---

[34] We note that, by definition, every recursive number-theoretic relation has denumerable formal representations in P (cf. [Me64], p118); for instance, if [$F$] is one such representation that satisfies Gödel's Theorem V ([Go31a], p22), so is [$F$ & (0=0)].

[35] We generally follow Mendelson's terminology and definitions of the "expressibility" ([Me64], §2, p117), and "representability" ([Me64], §2, p118), of number-theoretic relations and functions, respectively, in a formal system such as P. However, following Gödel, we also refer to a number-theoretic relation as being "representable" in P when, strictly speaking, we mean that it is "expressible" in P as defined by Mendelson.

[36] Loosely speaking, we cannot give a set-theoretic definition, of $Sb(x\ v|Z(y))$, such that $\{x\ |\ x=Sb(y\ 19|Z(y))\}$ defines a set in any Axiomatic Set Theory with a Replacement Axiom, or its equivalent, without introducing inconsistency.

A significant consequence is that we cannot consistently assume every recursive function formally defines a recursively enumerable set. It follows that we are unable to define a recursive set as a recursively enumerable set whose complement is also recursively enumerable. In some cases, there may be no such complement.

It further follows that, if $F(y)$ is an arithmetical function such that $F(k) = Sb(k\ 19|Z(k))$ for any given $k$, the assertion that the expression $\{x\ |\ x=F(y)\}$ defines a formal set by the Replacement Axiom may require additional qualification.

[37] In other words, even if instantiationally equivalent, a number-theoretic function, and a set-theoretic function (which is defined as a mapping), may be characteristics of different mathematical concepts.

In Meta-lemma 2, we argue that, even if a primitive recursive relation, and the standard interpretation of its formal representation, are always equivalent in their instantiations, they are not always formally equivalent[38].

In §II-1, §II-2 and §II-3, we consider the immediate implications of Meta-theorem 1 for the concepts of intuitive and formal consistency that are based on classical interpretations of Gödel's reasoning. In Meta-theorem 2, we consider, and analyse, an argument to the effect that no P-formula asserts, under the standard interpretation, that P is consistent.

In §II-4, Meta-lemma 3, we conclude that every recursive function is not strongly representable in P; and in Meta-lemma 4 that, even if a recursive relation is equivalent to some arithmetical relation in each of its instantiations, such equivalence cannot always be formulated within a formal system of Arithmetic such as P (or standard PA).

In §II-5, we address a removable ambiguity in Tarski's classical definitions of "satisfiability" and "truth". We then consider how these definitions can be made constructive, and intuitionistically unobjectionable, by formulating a constructive expression of Church's Thesis.

In §II-6, we introduce some constructive definitions of classical concepts, and offer a constructive definition of uncomputable number-theoretic functions.

In §II-7, we define a neo-classical Turing machine NT as a natural, and constructive, extension of a classical Turing machine T. We then introduce definitions for self-terminating, converging and oscillating NT-routines, and show that the significance of explicitly defining the truth of a formula of P under interpretation in terms of terminating

---

[38] We define two number-theoretic relations, say $f(x)$ and $g(x)$, as formally equivalent in P if, and only if, the equivalence "$f(x) <=> g(x)$" is the standard interpretation of some P-formula.

routines, of expressing Church's Thesis constructively, and of defining self-terminating computations of an NT machine, is expressed by the following meta-lemma:

Meta-lemma 14: If we assume a constructive Uniform Church Thesis, then every partial recursive number-theoretic function $F(x_1, ..., x_n)$ has a unique constructive extension as a total function.

We conclude that a constructive Uniform Church Thesis implies, firstly, that the classical Halting problem is effectively solvable (Corollary 14.1); and, secondly, that the classical Turing Thesis is false, and the Church-Turing equivalence is based on an invalid argument (Corollary 14.2 and Corollary 14.3, respectively).

In §II-8, Meta-lemma 17, and §II-9, Meta-lemma 18, we then show, respectively, that Turing's and Cantor's diagonal constructions do not necessarily determine Cauchy sequences[39].

## 2. Notation

Unless specified otherwise, we generally follow the notation introduced by Mendelson in his English translation of Gödel's 1931 paper [Go31a]; however, for convenience of exposition, we refer to it as Gödel's notation. Three notable exceptions: we use the notation "$(Ax)$", whose classical, standard, interpretation is "for all $x$", to denote Gödel's special symbol for Generalisation; the successor symbol is denoted by "$S$", instead of by "$f$"; and we use the symbol "¶" to denote the end of a proof.

Following Gödel (cf. [Go31a], footnote 13), we use square brackets to indicate that the expression $[(Ax)]$, including square brackets, only denotes the uninterpreted string[40] named[41] within the brackets. Thus, [(A$x$)] is not part of the formal system P, and would be replaced by Gödel's special symbolism for Generalisation wherever it occurs.

---

[39] Hence we cannot assume, as Turing does in the original formulation of his Halting problem ([Tu36], §8), that every "circle-free" Turing-machine necessarily defines a Dedekind real number. A similar objection holds for the classical assumption that Cantor's diagonal argument necessarily defines a Dedekind real number.

Following Gödel's definitions of well-formed formulas[42], we note that juxtaposing the string [(A$x$)] and the formula[43] [$F(x)$] is the formula [(A$x$)$F(x)$], juxtaposing the symbol [~] and the formula [$F$] is the formula [~$F$], and juxtaposing the symbol [v] between the formulas [$F$] and [$G$] is the formula [$F$v$G$].

The number-theoretic functions and relations in the following are defined explicitly by Gödel [Go31a]. The formulas are defined implicitly by his reasoning.

## 3. Definitions

We take P to be Gödel's formal system[44], and define ([Go31a], Theorem VI, p24-25):

(*i*)    $Q(x, y)$ as Gödel's recursive number-theoretical relation ~$xB(Sb(y\ 19|Z(y)))$[45].

(*ii*)   [$R(x, y)$] as a formula that represents $Q(x, y)$ in the formal system P.

---

[40] We define a "string" as any concatenation of a finite set of the primitive symbols of the formal system under consideration.

[41] We note that the "name" inside the square brackets only serves as an abbreviation for some string in P.

[42] We note that all well-formed formulas of P are strings of P, but all strings of P are not well-formed formulas of P.

[43] By "formula", we shall henceforth mean a "well-formed formula" as defined by Gödel ([Go31a], p11).

[44] Gödel ([Go31a], p9-13).

[45] We follow Gödel's definition of recursive number-theoretic functions and relations ([Go31a], p14-17). We note, in particular, that Gödel's recursive number-theoretic function $Sb(x\ 19|Z(y))$ is defined as the Gödel-number of the P-formula that is obtained from the P-formula whose Gödel-number is $x$ by substituting the numeral [$y$], whose Gödel-number is $Z(y)$, for the variable whose Gödel-number is 19 wherever the latter occurs free in the P-formula whose Gödel-number is $x$ ([Go31a], p20, def.31). We also note that Gödel's recursive number-theoretic relation $xBy$ holds if, and only if, $x$ is the Gödel-number of a proof sequence for the P-formula whose Gödel-number is $y$ ([Go31a], p22, def. 45).

(*iii*)  *q* as the Gödel-number[46] of the formula [*R*(*x*, *y*)] of P.

(*iv*)  *p* as the Gödel-number of the formula [(A*x*)*R*(*x*, *y*)][47] of P.

(*v*)  [*p*] as the numeral that represents the natural number *p* in P.

(*vi*)  *r* as the Gödel-number of the formula [*R*(*x*, *p*)] of P.

(*vii*)  17*Genr* as the Gödel-number of the formula [(A*x*)*R*(*x*, *p*)] of P.

(*viii*)  *Neg*(17*Genr*)[48] as the Gödel-number of the formula [~(A*x*)*R*(*x*, *p*)] of P.

(*ix*)  *R*(*x*, *y*) as the standard interpretation of the formula [*R*(*x*, *y*)] of P.

(*x*)  *Wid*(P) as the number-theoretic assertion (E*x*)(*Form*(*x*) & ~*Bew*(*x*))[49].

(We note that *Wid*(P) is defined by Gödel ([Go31a], p36) as equivalent to the meta-assertion "P is consistent".)

(*xi*)  [*Con*(P)] as the formula that represents *Wid*(P) in the formal system P.

(*xii*)  *w* as the Gödel-number of the formula [*Con*(P)] of P [1, p37].

(*xiii*)  *Con*(P) as the standard interpretation of the formula [*Con*(P)] of P.

---

[46] By the "Gödel-number" of a formula of P, we mean the natural number corresponding to the formula in the 1-1 correspondence defined by Gödel ([Go31a], p13).

[47] We note that "[(A*x*)][*R*(*x*, *y*)]" and "[(A*x*)*R*(*x*, *y*)]" denote the same formula of P.

[48] We note that Gödel's recursive number-theoretic function *Neg*(*x*) is the Gödel-number of the P-formula that is the negation of the P-formula whose Gödel-number is *x* ([Go31a], p18, def. 13).

[49] We note that Gödel's recursive number-theoretic relation *Form*(*x*) is satisfied if, and only if, *x* is the Gödel-number of a P-formula ([Go31a], p19, def. 23). Also, Gödel's number-theoretic relation *Bew*(*x*) is satisfied if, and only if, *x* is the Gödel-number of a provable P-formula ([Go31a], p22, def. 46).

## 4. A Meta-theorem of recursive asymmetry

**Definition 4**(*i*): A recursive number-theoretic function or relation is *asymmetrical* in P if it is not the standard interpretation of any of its formal representations in P.

**Meta-theorem 1**: There is a recursive number-theoretic relation that is asymmetrical in P.

*Proof*: We consider the number-theoretic relation $x=Sb(y\ 19|Z(y))$[50].

(*a*) We assume that no recursive number-theoretic function or relation is asymmetrical in P. In other words, we assume that every recursive number-theoretic function or relation is the standard interpretation of at least one of its formal representations in P.

(*b*) Let the P-formula $[F(x, y)]$ denote a formal representation of the recursive number-theoretic relation $x=Sb(y\ 19|Z(y))$.

(*c*) Let $F(x, y)$ denote the standard interpretation of $[F(x, y)]$.

(*d*) We consider the case where $x=Sb(y\ 19|Z(y))$ is an abbreviation[51] for $F(x, y)$.[52]

---

[50] This relation occurs on 2nd August 2002 in the correspondence titled "The Godel's Loop" between Antonio Espejo <CASAFARFARA@terra.es> and Rupert<rupertmccallum@yahoo.com> in the Google newsgroup sci.logic.

[51] Mendelson ([Me64], p31, footnote 1).

[52] In other words, we assume that, if we use Gödel's recursive definitions ([Go31a], p17-20), and follow the reasoning he outlines in Theorem V ([Go31a], p23), we can transform the relation $x=Sb(y\ 19|Z(y))$ into a relation $F(x, y)$, such that all the symbols that occur in $F(x, y)$ are standard interpretations of primitive symbols of P. See also Gödel's remarks ([Go31a], p11, footnote 22) in this context.

Now:

(*i*)  Let $k$ be the Gödel-number of [$F(x, y)$].

(*ii*)  Then $Sb(k\ 19|Z(k))$ is the Gödel-number of the P-formula [$F(x, k)$] that we get when, in the P-formula [$F(x, y)$], we replace the variable [$y$], wherever it occurs, by the numeral [$k$].

(We assume that [$y$] is Gödel-numbered as $19^{53}$ in the Gödel-numbering that yields $k$ as the Gödel-number of [$F(x, y)$])

(*iii*)  Let $l = Sb(k\ 19|Z(k))$.

(*iv*)  We now note that, by definition, the unevaluated numerical expression obtained by substituting $k$ for $y$ in the number-theoretic function $Sb(y\ 19|Z(y))$, which we abbreviate as "$Sb(k\ 19|Z(k))$" [54], must contain an explicit bound $k'$ that is equal to, or larger than, $l$.

This follows from the constructive definition of Gödel's recursive functions ([Go31a], p17, footnote 34), and his Theorem IV ([Go31a], p16). Thus $Sb(y\ 19|Z(y))$ is of the form "$(ex)((x = < f_1(y))\ \&\ g_1(x, y)))$"; where $f_1(y)$ and $g_1(x, y)$ are a recursive function and a recursive relation, respectively, both of lower rank than that of $Sb(y\ 19|Z(y))$. Similarly, $f_1(y)$ is of the form "$(ex)((x = < f_2(y))\ \&\ g_2(x, y))$", etc.

---

[53] Following the convention set by Gödel in Theorems V and VI (cf. [Go31a], p22, footnote 38), we assign 17 as the Gödel-number of "$x$", and 19 as the Gödel-number of "$y$".

[54] We use inverted commas to denote that "$Sb(k\ 19|Z(k))$" refers to the unevaluated numerical expression that is obtained from the number-theoretic function $Sb(y\ 19|Z(y))$ by substituting the natural number $k$ for $y$ in $Sb(y\ 19|Z(y))$, and eliminating any abbreviations, so that the resulting expression only contains standard interpretations of the primitive symbols of P, as required by our hypothesis.

(Here, "$(ex)R(x)$" denotes the smallest natural number for which the relation $R(x)$ is satisfied, and "$=<$" denotes the relation of "equal to or less than".)

We thus have a finite sequence $f_1(y), f_2(y), ..., f_n(y)$ of recursive functions of decreasing rank, such that:

$$Sb(y\ 19|Z(y)) =< f_1(y) =< f_2(y) =< ... =< f_n(y),$$

where $n$ is less than or equal to the rank $r$ of $Sb(y\ 19|Z(y))$, and where $f_n(y)$ is of rank one.[55]

It follows that $f_n(y)$ occurs, as a bound, in the unabbreviated number-theoretic function whose abbreviation is $Sb(y\ 19|Z(y))$. By definition, $f_n(y)$ is, therefore, either a constant, or of the form "$y+q$" ([Go31a], p23), where $q$ is a natural number that depends on the ranks of $Sb(y\ 19|Z(y))$ and its defining functions.

(We note that the representation of $f_n(y)$ in P, as envisaged in Theorem V, would be Gödel's term of the first type ([Go31a], p10), denoted by $[S^q y]$, where '$S^q$' denotes the pre-fixing of the primitive (*successor*) symbol $[S]$ of P, to $[y]$, $q$ times.)

---

[55] Specifically, the argument here is essentially that, using the symbol # as abbreviation for an unspecified sequence of arguments, each of which is either the free variable $y$ or a recursive function with $y$ as its only free variable, and the symbol * as abbreviation for an unspecified recursive function, also with $y$ as its only free variable:

$f_1(y) =\ Sb_*(\#)$ ... ([Go31a], p20, def. 30)

$f_2(y) =\ Su(\#)$ ... ([Go31a], p20, def. 27)

$f_3(y) =\ Pr(\#)^*$ ... ([Go31a], p20, def. 5)

$f_4(y) =\ \#!^*$ ... ([Go31a], p20, def. 4)

$f_5(y) =\ S^*0$ ... since $\#! = S^*0$, where "$S$" denotes the successor symbol.

We thus have that:

$$l = Sb(k \ 19|Z(k)) = < f_n(k) = k'.$$

(*v*) Since the standard interpretation of [*F*(*x*, *k*)] is "*F*(*x*, *k*)", it follows, from (*d*), that "*F*(*x*, *k*)" is the number-theoretic relation whose abbreviation is "*x*= *Sb*($k$ 19|*Z*(*k*))".

(*vi*) By (*iv*), [*F*(*x*, $k$)] must, therefore, contain a numeral [$k'$], which interprets as a natural number $k'$ that is larger than $l$.

This is impossible, since a formula cannot contain a numeral that, under interpretation, yields a natural number that is equal to, or larger than, the Gödel-number of the formula[56]. It follows that assumption (*a*) does not hold; this proves the meta-theorem.¶

## 5. Two meta-lemmas

It now follows that:

**Meta-lemma 1**: We cannot introduce a finite number of arbitrary recursive number-theoretic functions and relations, as function letters and predicate letters respectively, into P without risking inconsistency.

*Proof*: Adding "Sb", and a finite number of other functions and relations in terms of which it is defined ([Go31a], p17-19, def. 1-31), as new function letters and

---

[56] This can be easily proved since, following, for instance, Gödel's assignment of the natural numbers 1 and 3 to the primitive P-symbols [0] and [*S*] respectively ([Go31a], p13), we have that:

(*i*) the Gödel number $p_1^3 p_2^3 \ldots p_q^3 p_{q+1}$ of the numeral [$S^q 0$], which represents the natural number $q$ in P, is greater than $q$ for all $q >= 0$, where $p_i$ is the $i$'th prime, and

(*ii*) for any P-formulas [*F*] and [*G*], the Gödel number of the concatenated P-formula [*FG*] is always greater than, or equal to, the individual Gödel numbers of the P-formulas [*F*] and [*G*].

predicate letters, respectively, to P, along with associated defining axioms (cf.
[Me64], §9, p82), would yield a formal system P' in which [$x=Sb(y$ 19|$Z(y))$] is a P'-
formula. By Meta-theorem 1, P' would be inconsistent.¶

**Corollary 1.1**: There is a recursive number-theoretic function whose range does not
effectively well-define a (recursively enumerable ([Me64], p250)) sub-set of the set
of natural numbers in any Axiomatic Set Theory Z (cf. [Me64], p192-4) that contains
an Axiom of Replacement[57], and which is a model of P.

*Proof*: Let $f(y)$ be a recursive number-theoretic function, and [$F(y)$] its
representation in P. We define the number-theoretic relations $t(x)$ by:

$t(x) <=> (Ey)(x=f(y))$.

Since a number-theoretic relation is expressible in P if, and only if, it is recursive
([Me64], Corollary 3.29, p142), and it can be shown that $t(x)$ is not always recursive
([Me64], Proposition 5.20(5), p251), it follows that $t(x)$ is not always expressible in
P.

Assuming that $t(x)$ is not expressible in P, we note that the set $T$, defined as $\{x \mid t(x)\}$, is, nevertheless, a well-defined sub-set, of the set of natural numbers, which is
definable in Z. This follows since:

(*i*)   for any given natural number $k$, the standard, arithmetical, interpretation
        $F(y)$ of the P-formula [$F(y)$] is such that, by definition:

        $f(k) = F(k)$);

---

[57] We take the Replacement Axiom of an Axiomatic Set Theory Z as stating, essentially, that the range of
every function of Z, whose domain is a well-defined set in Z, and whose values are always elements of a
well-defined set in Z, is a well-defined set in Z.

(*ii*)  the set *T′*, defined by {*x* | *t′*(*x*)}, is a well-defined set in Z, where:

$$t'(x) <=> (\mathrm{E}y)(x=F(y));$$

(*iii*)  by the Axiom of Replacement and (*i*), {*x* | *t*(*x*)} is a sub-set of {*x* | *t′*(*x*)}.

It follows that we can add a new function letter *f* to Z, without inviting inconsistency, such that:

*x* is in {*x* | *t*(*x*)} if, and only if, (E*y*)(*x*=*f*(*y*)).

Now, since Z is a model for P, it further follows that we can also add the function letter *f* to P, along with suitable defining axioms, without inviting inconsistency. However, if we take *f*(*y*) as the recursive number-theoretic function *Sb*(*y* 19|*Z*(*y*)), it follows by Meta-lemma 1 that such addition would introduce inconsistency into P.

We conclude that the range of the recursive number-theoretic function *Sb*(*y* 19|*Z*(*y*)) does not effectively well-define a (recursively enumerable) sub-set of the set of natural numbers in Z.¶

**Corollary 1.2**: Not every constructively well-defined number-theoretic function is a mathematical object.

It further follows:

**Meta-lemma 2**: Even if a primitive recursive relation, and the standard interpretation of its formal representation, are always equivalent in their instantiations, they are not always formally equivalent[58].

---

[58] We define two number-theoretic relations, say *f*(*x*) and *g*(*x*), as formally equivalent in P if, and only if, the equivalence "*f*(*x*) <=> *g*(*x*)" is the standard interpretation of some P-formula.

*Proof*: Let $F(x_1, ..., x_n)$ be any recursive number-theoretic relation, and $G(x_1, ..., x_n)$ be the standard interpretation of one of its formal representations in P.

We assume that the two relations are equivalent in their instantiations, so that, for any given sequence $<a_1, ..., a_n>$ of natural numbers:

$F(a_1, ..., a_n)$ holds if, and only if, $G(a_1, ..., a_n)$ holds.

However, by Meta-lemma 1, it follows that $[F(x_1, ..., x_n)]$ is not necessarily a P-formula. Hence, we cannot conclude that:

$[F(x_1, ..., x_n) <=> G(x_1, ..., x_n)]$ is a P-formula.[59]¶

## II. Consequences

### 1. Analysing Gödel's Theorem XI

Now we note that Gödel's number-theoretic relation *Bew*(*x*) is defined ([Me64], p22) in terms of his recursive number-theoretical relation *Sb*(*x* 19|*Z*(*y*)). The question thus arises:

If we assume that the number-theoretic sentence (E*x*)(*Form*(*x*) & ~*Bew*(*x*)) ([Me64], p36, footnote 63), abbreviated *Wid*(P), defines the proposition "P is consistent" in a constructive and intuitionistically unobjectionable way[60], then can we consistently assume further that *Wid*(P) is equivalent to the standard interpretation of some P-formula [*Con*(*P*)]?

---

[59] We note that Gödel's Theorems VIII to XI are based on the premise that the equivalence in *Meta-lemma 2* can always be formulated within the formal system P ([Go31a], p31).

[60] In other words, this definition can be assumed equivalent to Mendelson's classical meta-definition of consistency ([Me64], p37). We argue in §II-2 that this assumption, too, may conceal an implicit ambiguity.

We note that the latter assumption is one of the implicit meta-theses that appear to underlie Gödel's proof of, and the conclusions he draws from, his Theorem XI ([Go31a], p36).

Clearly, the reasoning in Meta-theorem 1 and Meta-lemma 2 implies that such an assumption is invalid[61]. We conclude that there is no P-formula, Con(P), whose standard interpretation is the number-theoretic assertion (E$x$)($Form(x)$ & $\sim Bew(x)$) - which is defined by Gödel as equivalent to "P is consistent".

## 1.1 Implicit meta-theses underlying Gödel's Theorem XI

The question then arises: Is there any P-formula whose standard interpretation can be defined equivalent to the proposition "P is consistent" in a constructive, and intuitionistically unobjectionable, way?

In Meta-theorem 2 we address this question by considering an argument that is based on the assumption that P formalises Dedekind's Peano Arithmetic[62]. However, we do not appeal to the further thesis that P is consistent if it has a model, since the consistency of the model may appeal to set-theoretic reasoning that is non-constructive, and intuitionistically objectionable. Instead, without assuming the consistency of P, we appeal to the following four meta-theses, all of which appear to be implicit in, or consequences of, Gödel's reasoning in the proof of his Theorem XI ([Go31a], p36):

---

[61] For, if (E$x$)($Form(x)$ & $\sim Bew(x)$) is formally equivalent to the standard interpretation of one of its formal representations, then so also are each of $Form(x)$ and $Bew(x)$. Arguing similarly, this would eventually imply that the recursive function $Sb(x\ 19|Z(y))$ too is the standard interpretation of one of its formal representations, contradicting Meta-theorem 1.

[62] This assumption is stated as an explicit, albeit informal, meta-premise by Gödel in his 1931 paper ([Go31a], p10): "P is essentially the system which one obtains by building the logic of PM around Peano's axioms ...".

**Meta-thesis 1**: P is a faithful formalisation of Dedekind's Peano Axioms[63]

**Meta-thesis 2**: If a sentence [*F*] is P-provable, then its standard interpretation *F* is a true number-theoretic sentence.[64]

**Meta-thesis 3**: "P is consistent", abbreviated *Wid*(P), can be defined as some number-theoretic sentence in a constructive, and intuitionistically unobjectionable, way.

**Meta-thesis 4**: *Wid*(P) is equivalent to the standard interpretation, *Con*(P), of some P-formula [*Con(P)*].

## 1.2 A negative meta-theorem

We now argue that:

**Meta-theorem 2**: If we assume that Meta-thesis 1-3 are true, then no P-formula can assert P as consistent under the standard interpretation.

*Proof*: We assume that Meta-theses 1-4 are true, so that there is some formula [*Con*(P)] of the formal system P such that, under the standard interpretation:

---

[63] We note that, in his 1931 paper [Go31a], Gödel explicitly assumes that P formalises "... the ordinary methods of definition and proof of classical mathematics..." ([Go31a], p36), and that every proof sequence, of P, interprets as a finite sequence of true assertions of Dedekind's Peano Arithmetic. Thus every member of the sequence is either true, since it is the interpretation of an axiom, or it is true, since it follows by classical logic from the previous true assertions in the sequence. Hence, every provable formula of P interprets as a true assertion under the standard interpretation.

[64] Meta-thesis 2 is clearly a consequence of the assumption that P is consistent. However, as we argue in §II-5, if we follow Tarski's definitions of the "satisfiability" and "truth" of provable sentences of P under an interpretation ([Me64], p50-53), then we cannot assume that Dedekind's Peano Arithmetic is consistent in a constructive, and intuitionistically unobjectionable way (cf. §II-5(b)). Hence, the converse need not be true; in other words, we cannot constructively assume Meta-thesis 2 implies P is consistent. Prima facie, Meta-thesis 2 differs from the assertion that P is classically sound; it does not assume, or imply, that a sentence can be true only under an interpretation that is intuitively consistent.

Con(P) holds if, and only if, Wid(P) holds.

We take this as equivalent, by Meta-thesis 3, to the assertion:

[Con(P)] is a P-formula that asserts, under the standard interpretation, that P is consistent.

(*i*)   By the classical definition of consistency, [Con(P)] is P-provable if P is inconsistent - since every formula of an inconsistent P is a consequence of the Axioms, by the Rules of Inference, of P[65].

(*ii*)  Now, if [Con(P)] were P-provable then, by Meta-thesis 2, we would conclude, under the standard interpretation, that:

Con(P) is a true number-theoretic sentence.

We would further conclude, by Meta-thesis 3 and Meta-thesis 4, the meta-assertion:

P is a consistent formal system.

However, this would be an invalid meta-conclusion, since P may be inconsistent.

(*iii*) It follows that we cannot conclude from the P-provability of [Con(P)] that P is consistent[66]. Hence we cannot have any formula [Con(P)] such that, under the standard interpretation:

If Con(P) is a true number-theoretic sentence, then Wid(P) is a true number-theoretic sentence.

---

[65] This follows from ([Me64], p37, Corollary 1.15).

[66] As we note in the Analytic summary below, this argument may be only apparently circular.

(*iv*)  We conclude that, if Meta-thesis 1-3 are assumed true, then Meta-thesis 4 is false, and there is no P-formula [*Con*(P)] such that, under the standard interpretation:

*Con*(P) is a true number-theoretic sentence if, and only if, P is a consistent formal system.¶

The question arises: How does Meta-theorem 2 affect Gödel's Theorem XI ([Go31a], p36)?

## 1.3  Gödel's Proof of Theorem XI

Now Gödel's Theorem XI is, essentially, the following proposition.

**Gödel's Theorem XI**: The consistency of P is not provable in P if P is consistent.

*Proof*:  Gödel argues that:

(*i*)  If P is assumed consistent, then the following number-theoretic assertions follow from his Theorems V, VI and his definition of *Wid*(P).

$Wid(P) \Rightarrow \sim Bew(17Genr)$

$Wid(P) \Rightarrow (Ax)\sim xB(17Genr)$

$17Genr = Sb(p\ 19|Z(p)))$

$Wid(P) \Rightarrow (Ax)\sim xB(Sb(p\ 19|Z(p)))$

$Q(x, y) <=> \sim xB(Sb(y\ 19|Z(y)))$

$(x)Q(x, y) <=> (Ax)\sim xB(Sb(y\ 19|Z(y)))$

$Wid(P) \Rightarrow (Ax)Q(x, p)$

(*ii*)  Assuming that [(A*x*)*R*(*x, p*)] asserts its own provability, Gödel concludes from the above that the instantiation:

  *wImp*(17*Genr*),

of the recursive number-theoretic relation *xImpy*, is a true number-theoretic assertion under the standard interpretation.

(*iii*)  From this, he concludes that:

  [*Con*(P) => (A*x*)*R*(*x, p*)] is P-provable.

(*iv*)  Now, in his Theorem VI, Gödel (1) argues that, if P is assumed consistent, then [(A*x*)*R*(*x, p*)] is not P-provable. He thus concludes that, if P is assumed consistent, then [*Con*(P)] too is not P-provable.

(*v*)  Implicitly assuming that Meta-thesis 4 is true, and so, under the standard interpretation:

  *Con*(P) is a true number-theoretic sentence if, and only if, *Wid*(P) is a true number-theoretic sentence,

Gödel further concludes that (*iv*) is equivalent to asserting that the classical notion of the consistency of P is not provable in P.¶

## 1.4  Gödel's Theorem XI as a conditional meta-assertion

However, since Gödel's Theorem XI implicitly assumes Meta-theses 1-4, it should be treated, more appropriately, as the conditional meta-assertion:

If, assuming Meta-theses 1-4[67], there is a P-formula [*Con*(P)] whose standard interpretation is equivalent to the assertion "P is consistent", then [*Con*(P)] is not P-provable if P is consistent.

Hence Meta-theorem 2 implies that Gödel's Theorem XI is a vacuous meta-assertion.

## 1.5  Analytical summary

(*a*) What we have essentially argued above is that, for his Theorem XI to hold, Gödel needs to, firstly, meta-establish both (*i*) and (*ii*):

(*i*)   P is a consistent formal system

==>[68] "(E$x$)(*Form*($x$) & ~*Bew*($x$))" is a true[69] number-theoretic assertion

==> "*Con*(P)" is a true number-theoretic assertion

(*ii*)  "*Con*(P)" is a true number-theoretic assertion

==> "(E$x$)(*Form($x$)* & ~*Bew*($x$))" is a true number-theoretic assertion

==> P is a consistent formal system

Now what Gödel actually meta-argues is (*iii*) and (*iv*):

---

[67] We note, again, that "P is consistent" implies Meta-thesis 2.

[68] We use the long-arrow notations "==>" and "<==>" as abbreviations for "implication" and "equivalence", respectively, between meta-propositions.

[69] Where it is obvious from the context, we assume that "true" means "true under the standard interpretation".

(*iii*) P is a consistent formal system

==> "(E*x*)(*Form*(*x*) & ~*Bew*(*x*))" is a true number-theoretic assertion ... (*by definition*)

==> [*Con*(P)] is not P-provable.

(*iv*) "(E*x*)(*Form*(*x*) & ~*Bew*(*x*))" is a true number-theoretic assertion

==> P is a consistent formal system ... (*by definition*)

From this he implicitly concludes, in Theorem XI, that:

(*v*) "P is consistent" is equivalent to [*Con*(P)] under the standard interpretation, and that [*Con*(P)] is not P-provable.

However, since he does not meta-establish that:

(*vi*) P is a consistent formal system

<==> "*Con*(P)" is a true number-theoretic assertion

his conclusion in Theorem XI can only be that:

(*vii*) If "P is consistent" is equivalent to [*Con*(P)] under the standard interpretation, then [*Con*(P)] is not P-provable.

(*b*) In Meta-theorem 2 above, we argue, as below, that (*a*)(*vii*) holds vacuously:

(*i*) Assuming that P is a faithful formalisation of Dedekind's Peano Axioms, the axioms of P are true under the standard interpretation by the classical definitions of "satisfiability" and "truth" for predicate logic (cf. [Me64], p50), and the rules of inference preserve truth. Thus, every provable formula of P is also true under the standard interpretation. It follows that:

[*Con*(P)] is P-provable ==> "*Con*(P)" is a true number-theoretic assertion.

(*ii*)  Now we have that:

P is an inconsistent formal system ==> [*Con*(P)] is P-provable

(*iii*)  If we, then, assume that:

"*Con*(P)" is a true number-theoretic assertion ==> P is a consistent formal system

it would follow from (*i*) that:

[*Con*(P)] is P-provable ==> P is a consistent formal system.

(*iv*)  This is clearly false, since, by (*i*) and (*ii*), we may have that:

[*Con*(P)] is P-provable & P is an inconsistent formal system[70].

(*v*)  We conclude that, assuming P is a faithful formalisation of Dedekind's Peano Axioms, there is no number-theoretic assertion "*Con*(P)" such that:

"*Con*(P)" is a true number-theoretic assertion ==> P is a consistent formal system.

---

[70] This, of course, would imply that Dedekind's Peano Arithmetic, too, is inconsistent. Assuming that such could be the case, it would follow that the inconsistency is prevented from becoming intuitively evident by some, yet more fundamental, implicit premises that may be involved in the semantic interpretation of the sentences of the Arithmetic.

## 2. Consistency and Meta-thesis 3

Prima facie, it may seem reasonable to accept Gödel's assertion - that "P is consistent" can be defined as equivalent to the number-theoretic sentence (E$x$)(*Form*($x$) & ~*Bew*($x$)) - and to assume that this definition is equivalent to the classical meta-definition of consistency, such as that expressed in Mendelson ([Me64], p57, Corollary 1.15).

However, the number-theoretic assertion (E$x$)(*Form*($x$) & ~*Bew*($x$)) can be treated as equivalent to the meta-assertion "P is consistent", only if we understand it non-constructively to mean:

> There is some natural number $n$ for which we can assert that $n$ is the Gödel -number of some formula [$F$] of P that is not P-provable.

Now, since "(E$x$)" is only a shorthand notation for "~(A$x$)~" ([Go31a], p11), the assertion "(E$x$)(*Form*($x$) & ~*Bew*($x$))" should actually be read as "~(A$x$)~(*Form*($x$) & ~*Bew*($x$))". As we argue in §II-5, we can interpret the above constructively, either as the meta-assertion:

> (*i*)  It is not true that, for any given natural number $n$, there is an *individually* constructive, and intuitionistically unobjectionable, way to determine that, it is not true that, $n$ is the Gödel-number of some formula [$F$] of P that is not P-provable.

or, as the meta-assertion:

> (*ii*) It is not true that there is a *uniformly* effective, and intuitionistically unobjectionable, way to determine that, for any given natural number $n$, it is not true that, $n$ is the Gödel-number of some formula [$F$] of P that is not P-provable.

We note that, although (*ii*) obviously implies (*i*), the converse need not be true. Thus, we cannot assume that the two meta-assertions are necessarily equivalent. So Gödel's assumption that the number-theoretic sentence (E*x*)(*Form*(*x*) & ~*Bew*(*x*)) is a precise, and intuitionistically unobjectionable, definition of "P is consistent" is not unambiguous; it can, conceivably, lead to anomalous consequences.[71]

We note that a similar ambiguity may exist in the interpretation of the consequences of Gödel's construction of the recursive relation $Q(x, p)$[72]; thus, although we can effectively assert that "$Q(n, p)$ holds individually for any given natural number *n*", we cannot, prima facie, assume that this is equivalent to the non-constructive, infinite, compound, sentence "$Q(x, p)$ holds uniformly for all natural numbers *x*".

The distinction between the two assertions is better expressed in terms of classical, program-terminating, Turing routines[73]. Thus, given any *n*, it follows - from Gödel's reasoning that $[R(n, p)]$[74] is P-provable - that there is always some effective method that will terminate in a finite, even if indeterminate, number of steps $t(n)$ if, and only if, $R(n, p)$ holds. However, since $[(Ax)R(x, p)]$ is not P-provable, there may not be any Turing machine such that, given any *n*, it will halt in a determinate number of steps $t(n)$[75] if, and only if, $R(n, p)$ holds[76]. In other words, there may be no classical, program-terminating,

---

[71] Another instance, albeit unrelated, where the intuitive interpretation of a recursive number-theoretic function could lead to possibly anomalous consequences, is Mendelson's "remainder" function $rm(x, y)$ ([Me64], p122, Proposition 3.15(n)). If, following Mendelson, we interpret $rm(x, y)$ as yielding the remainder upon division of *y* by *x*, we have the anomalous interpretation that the remainder upon division of *y* by 0 is *y*, since it follows from his definition that, for any natural number *y*, $rm(0, y) = y$!

[72] Cf. Gödel's recursive relation $Q(x, y)$ ([Go31a], p24, Eqn. 8.1).

[73] We take Mendelson's exposition ([Me64], p228) as representative of the classical definition of "Turing algorithms"; we also refer to a classical "Turing algorithm" as a "program-terminating Turing routine".

[74] We note that $[R(x, p)]$ formally expresses the recursive relation $Q(x, p)$ in P.

[75] By a "step", we mean here an instantaneous tape description of a Turing machine ([Me64], p230).

[76] This would, for instance, be the case if, for any Turing machine T that halts on input *n* only if $R(n, p)$ holds, there is some input *n* on which T loops, but $R(n, p)$ holds.

Turing routine for computing the function $t(x)$; in Turing's terminology, $t(x)$ may be Turing-uncomputable, even though $t(n)$ is effectively computable for any given $n$.[77]

## 3. Can consistency be a formal convention?

Meta-theorem 2 and §II-2 raise the question: Can we arbitrarily postulate that a formal sentence, under an interpretation, asserts the consistency of the system? If so, the formula expressing such consistency may be in the nature of an impermanent convention that has intuitive significance only for a transient panel of mathematical logicians. That the issue is not trivial is indicated by Mendelson's remarks ([Me64], p148):

> "Let $Con$(S) be the wf: $(x_1)(x_2)(x_3)(x_4)\sim(Pf(x_1, x_3)\ \&\ Pf(x_2, x_4)\ \&\ Ng(x_3, x_4))$. Intuitively, according to the standard interpretation, $Con$(S) asserts that there is no proof in S of any wf and its negation ... ."[78]

Mendelson then goes on to remark ([Me64], p149), for a first order theory K that possesses the individual constants of S:

> "... the way in which $Con$(K) is constructed also adds an element of ambiguity. This ambiguity is dangerous, because ... there is a reasonable way of defining $Con$(S) so that $Con$(S) is S-provable."

Mendelson finally outlines Feferman's formulation of "... a wf expressing the consistency of K". However, intuitively, this formula no longer even appears related to Mendelson's classical, and intuitionistically unobjectionable, notion of consistency ([Me64], p57, Corollary 1.15).

---

[77] This argument is developed in Corollary 14.1.

[78] Even if we accept that, prima facie, this definition adequately formalises our notion of consistency, it implicitly assumes that the formulas $Pf(x, y)$ and $Ng(x, y)$ do, indeed, interpret as semantically - and not merely instantiationally - equivalent to the recursive relations that they represent. As we remark in the concluding paragraph of the previous section, this, itself, is a significant assumption of questionable validity..

The question arises: If these definitions cannot be intuitively interpreted as representing the classical expression of consistency, then is there a constructive, and intuitionistically unobjectionable, logical basis by which a formal postulation can be meaningfully adopted as an effective definition to replace the classical expression of consistency?[79]

## 4. Definitions of new function and predicate letters

Meta-theorem 1 highlights a non-constructive issue underlying Gödel's reasoning in his Theorem V. It addresses the question: Can every recursive number-theoretic function or relation be introduced effectively as a formal mathematical object into P?

### 4.1 The classical argument

The classical argument is expressed by Mendelson ([Me64], p82, §9):

"In mathematics, once we have proved, for any $y_1, ..., y_n$, the existence of a unique object $u$ having the property $A(u, y_1, ..., y_n)$, we often introduce a new function $f(y_1,$

_____________________

[79] We note that, referring generally to a somewhat similar issue concerning whether axioms of a formal theory need to be intuitively self-evident, Gödel argues in favor of adding to mathematics axioms which are not self-evident and which are only justified pragmatically: ". . . even disregarding the intrinsic necessity of some new axiom, and even in case it has no intrinsic necessity at all, a probable decision about its truth is possible also in another way, namely, inductively by studying its "success." Success here means fruitfulness in consequences, in particular in "verifiable" consequences, i.e., consequences demonstrable without the new axiom, whose proofs with the help of the new axiom, however, are considerably simpler and easier to discover, and make it possible to contract into one proof many different proofs. The axioms for the system of real numbers, rejected by intuitionists, have in this sense been verified to some extent, owing to the fact that analytical number theory frequently allows one to prove number-theoretical theorems which, in a more cumbersome way, can subsequently be verified by elementary methods. A much higher degree of verification than that, however, is conceivable. There might exist axioms so abundant in their verifiable consequences, shedding so much light upon a whole field, and yielding such powerful methods for solving problems (and even solving them constructively, as far as that is possible) that, no matter whether or not they are intrinsically necessary, they would have to be accepted at least in the same sense as any well-established physical theory." (Kurt Gödel, 1947, "What is Cantor's continuum problem?", Collected Works, vol. 2, Oxford University Press, 1986–2003.)

..., $y_n$) such that $A(f(y_1, ..., y_n), y_1, ..., y_n)$ holds for all $y_1, ..., y_n$. ... It is generally acknowledged that such definitions, though convenient, add nothing really new to the theory."

More precisely, he argues, in his Proposition 2.29 ([Me64], p82), that, classically, abbreviations for strongly representable ([Me64], p118) number-theoretic functions may be introduced as function letters into P, since they can always be eliminated as follows:

Let K be a first-order theory with equality. Assume that $[(E!u)A(u, y_1, ..., y_n)]$[80] is K-provable. Let K' be the first-order theory with equality obtained by adding to K a new function letter $[F]$ of $n$ arguments, and the proper axiom $[A(F(y_1, ..., y_n), y_1, ..., y_n)]$, as well as all instances of the axioms of K involving $[F]$. Then:

(*i*)   If K' is assumed consistent[81], every formula $[B]$ of K' is provably equivalent to a formula $[B']$ of K' that does not contain $[F]$.

(*ii*)   If K' is assumed consistent[82], $[B']$ is K-provable if, and only if, $[B']$ is K'-provable.

However, assuming that every recursive number-theoretic function can be strongly represented in P[83], it follows from Meta-theorem 1 that such classical definitions, as envisaged in Mendelson's Proposition 2.29, may introduce inconsistency into a consistent P. Hence K' cannot, then, be assumed consistent.

---

[80] We use the notation "E!", instead of Mendelson's "$E_1$", to denote "unique existence" ([Me64], p79).

[81] The assumption of consistency lies in the implicit use of argument by contradiction.

[82] The assumption of consistency lies in the implicit use of argument by contradiction.

[83] By Mendelson's reasoning in his Proposition 3.23 ([Me64], p131), every recursive function is representable in P. Mendelson also observes that it is not known whether every recursive function is strongly representable in P ([Me64], p135, Ex. 3).

We thus have:

**Meta-lemma 3**: If every strongly representable number-theoretic function can be introduced as a new function letter into P, without affecting the consistency of P, then every recursive function is not strongly representable in P.

## 4.2 Non-constructive definitions

Mendelson's remarks can be taken as reflecting an implicit, and non-constructive, classical belief, which his Proposition 2.29 apparently attempts to formalise explicitly. This is that, if we can represent every instantiation of a recursive function individually in K, it is sufficient to well-define the function as a formal, mathematical entity in K.

However, as noted in §II-2, in order to address this issue in a constructive, and intuitionistically unobjectionable way, classical arguments may need to explicitly address the distinction between:

(*i*) constructing, for any given natural number $n$, an individual decision-routine that terminates in a finite, even if indeterminate, number of steps on a given condition, and,

(*ii*) constructing a uniform decision-routine that will terminate in a determinate number of steps on a given condition for any given natural number $n$.

## 4.3 Individually and uniformly computable functions

To highlight the distinction, we assume that $[(E!y)R(y, x_1, ..., x_n)]$ is P-provable. The classical, non-constructive, conclusion drawn from this is that:

(*i*)   There is a unique $b$ in the domain of the interpretation M of P such that, given any sequence $<a_1, ..., a_n>$ in the domain of M, $R(b, a_1, ..., a_n)$ holds in M, where $R(y, x_1, ..., x_n)$ is the interpretation of $[R(y, x_1, ..., x_n)]$ in M.

We may thus define a total function $F(x_1, ..., x_n)$ in M such that, given any

sequence $<b, a_1, ..., a_n>$ in M:

$$R(b, a_1, ..., a_n) <=> (F(a_1, ..., a_n) = b)$$

holds in M.

Thus, (*i*) assures us, classically, that:

(*ii*)   For each given sequence $<a_1, ..., a_n>$ in the domain of M, there is some individual computing-routine, in M, that will compute a unique $b$ as the value of $F(a_1, ..., a_n)$, in a finite, possibly indeterminate, number of steps $t(a_1, ..., a_n)$.

However, the P-provability of $[(E!x)R(x, y_1, ..., y_n)]$, by itself, does not assure us that

there is a well-definable, uniformly terminating routine[84] in M, such that, given any

sequence $<a_1, ..., a_n>$ in the domain of M, it will compute a unique $b$ in the domain of

M, as the value of $F(a_1, ..., a_n)$, in a determinate number of steps $t(a_1, ..., a_n)$.

In other words, some sequence $<a_1, ..., a_n>$ may require a sequence-specific terminating routine in order to compute $F(a_1, ..., a_n)$. Hence, in the absence of a uniformly terminating routine, $F(x_1, ..., x_n)$ may be defined as individually, but not uniformly, computable in M.

---

[84] We treat the general concepts of a "terminating routine", a "terminating decision-routine", or a "terminating computing-routine", as essentially intuitive. Specifically, we would treat a "finite, terminating decision-routine" in P as equivalent to a proof sequence of P. A Cauchy sequence in a metric space could also be viewed intuitively as an "infinite, terminating computing-routine".

(There is also no assurance of a well-definable terminating routine in M that, for any given sequence $<a_1, ..., a_n>$, will compute $t(a_1, ..., a_n)$.)

## 4.4  Gödel's Theorem V

In view of the above, a reasonable interpretation of Gödel's outline of the proof of Theorem V ([Go31a], p22) indicates that, it, too, may implicitly rely on non-constructive reasoning in order to assert the formal existence of a unique mathematical object corresponding to any recursive number-theoretic function or relation.

Classically, such an assertion would imply that every recursive number-theoretic relation can either be expressed constructively, and in an intuitionistically unobjectionable manner, in terms of only symbols that are the interpretations of the primitive symbols of P, or introduced as a formal mathematical object into P without inviting inconsistency. Gödel's remarks, both whilst outlining a proof of Theorem V[85], and those at the end of his Theorem VII, indicate that this may have been his belief at the time. Such a belief would imply that the expression "$x=Sb(y\ 19|y)$" can be treated as an abbreviation of one of its representations in P, without affecting the consistency of P. Meta-theorem 1 argues that any such belief would be false.

## 4.5  Gödel's Theorem VII

We note that Gödel's Theorem VII ([Go31a], p29) asserts that, for any given set of natural number values of its free variables, every recursive number-theoretic relation, say $F(x)$, is equivalent to a corresponding instantiation of an arithmetical[86] relation $G(x)$ (i.e. a relation that uses only "+" and "*" as primitive arithmetic symbols). From this, Gödel

---

[85] "The definitional procedures by which f arises from $f_1$, ..., $f_k$ (substitution and recursive definition) can both be formally imitated in the system P" ([Go31a], p23).

[86] We use the term "arithmetical" as defined by Gödel ([Go31a], p29).

apparently concludes that the two relations must also be formally equivalent. He remarks that ([Go31a], p31):

> "According to Theorem VII, for every problem of the form $(x)F(x)$ ($F$ recursive), there is an equivalent arithmetical problem, and since the whole proof of Theorem VII can be formulated (*for each particular F*) within the system P. this equivalence is provable in P."

However, it follows, from Meta-lemma 1 and Meta-lemma 2, that such a premise is false. We cannot add Gödel's first 31 definitions to P - by adding suitable function and predicate letters, and treating their definitions as axioms[87] - without risking inconsistency. As such, we cannot assume that "the whole proof of Theorem VII can be formulated (*for each particular F*) within the system P". It follows that, where Gödel's reasoning in his Theorems VIII to XI appeals implicitly to the above premise, it may invite inconsistency.

## 4.6  Standard PA

We also note that every arithmetical relation is the standard interpretation of one of its representations in a formal system of Arithmetic such as standard PA[88]. Now, since Meta-Lemma 1 can be seen to hold in any such system, it follows from Meta-lemma 2 that:

> **Meta-lemma 4**:  Even if a primitive recursive relation is equivalent to some arithmetical relation, in the sense that they may be always equivalent in their instantiations, such equivalence cannot always be formulated within a formal system of Arithmetic such as standard PA.

---

[87] Cf. Mendelson ([Me64], p82, Proposition 2.29).

[88] We take standard PA to be the first  order theory S defined by Mendelson ([Me64], p102), in which addition and multiplication are the standard interpretations of the primitive symbols "+" and "*" respectively.

It also follows that Gödel's Theorems V to XI cannot be assumed proven, as he intended by his reasoning, in a formal system of Arithmetic such as standard PA.

## 5.  Classical definitions of "satisfiability" and "truth" are ambiguous

We next note that classical definitions, of "satisfiability" and "truth" ([Me64], p50-52), do not explicitly address the distinction between "*individually*[89] decidable" relations, and "*uniformly*[90] decidable" relations[91]. Thus, whereas a uniformly decidable relation is always individually decidable, Gödel's recursive number-theoretic relation $Q(x, y)$ may correspond to an instance where, under a constructive expression of Church's Thesis, the converse does not hold.

We consider removing such ambiguity by, firstly, introducing the additional terminology, "*terminating routine*", for the classical, intuitive, concept of an "effective method", or "mechanical procedure" ([Me64], p207).

## 5.1  Classical definitions of "satisfiability" and "truth"

We now note that the existence of a terminating routine, to effectively determine that the assertions involved in the following definitions hold in an interpretation M of the formal system P, is implicit in Tarski's classical definitions of the satisfiability and truth of the well-formed formulas of a formal language under a given interpretation ([Me64], p50-51). In order to highlight that these definitions may be verifiable by some effective method, we therefore introduce this condition explicitly as below, and consider some consequences.

---

[89] We use the terminology "individually decidable" if, and only if, a P-formula with free variables is decidable individually, as satisfiable or not, in M for any given instantiation. This corresponds to, but is not necessarily a consequence of, the case where every instantiation of the formula is P-provable. We note that the closure of the formula may not be P-provable.

[90] We use the terminology "uniformly" if, and only if, there is an effective method ([Me64], p207) by which a P-formula with free variables is decidable jointly, as true or not, in M as an infinitely compound assertion of all of its instantiations. This corresponds to, and is a consequence of, the case where the closure of the formula is P-provable.

[91] We note that Turing [Tu36] makes this distinction implicitly when he states that "Let P be a sequence whose $n$-th figure is 1 or 0 according as $n$ is or is not satisfactory. It is an immediate consequence of the theorem of §8 that P is not computable. It is (so far as we know at present) possible that any assigned number of figures of P can be calculated, but not by a uniform process. When sufficiently many figures of P have been calculated, an essentially new method is necessary in order to obtain more figures".

**Definition 5**(*i*): If $[R(x_1, ..., x_n)]$ is an atomic formula[92] of P, and $R(x_1, ..., x_n)$ is the corresponding relation of the interpretation, then *the sequence*[93] $<a_1, ..., a_n>$ of M *satisfies* $[R(x_1, ..., x_n)]$ *classically*, under the interpretation M, if, and only if, there is a terminating routine to determine that $R(a_1, ..., a_n)$ holds[94] in M, where $a_1, ..., a_n$ are elements in the domain D of M[95].

**Definition 5**(*ii*): If $[R]$ is a formula of P, a *sequence s of* M *satisfies* $[\sim R]$ *classically*, under the interpretation M, if, and only if, there is a terminating routine to determine that *s* does not satisfy $[R]$ classically.

**Definition 5**(*iii*): If $[R]$ and $[S]$ are formulas of P, *a sequence s of* M *satisfies* $[R{=}{>}S]$ *classically*, under the interpretation M, if, and only if, there is a terminating routine to determine that either *s* does not satisfy $[R]$ classically, or *s* satisfies $[S]$ classically.

**Definition 5**(*iv*): If $[R]$ is a formula of P, *a sequence s of* M *satisfies* $[(Ax_i)R]$ *classically*, under the interpretation M, if, and only if, there is a terminating routine to determine that every sequence of M, which differs from *s* in at most the *i*'th component, satisfies $[R]$ classically.

---

[92] We follow Mendelson's definition of an "atomic formula" ([Me64], p46).

[93] Unless otherwise specified, all sequences are assumed to be finite, ordered, sets of elements of the domain D of an interpretation M, and denoted as *<sequence>*.

[94] We take the meaning of "$R(a_1, ..., a_n)$ holds in M", where $[R(a_1, ..., a_n)]$ is an atomic formula, to be intuitively well-understood.

[95] We do not assume that the elements in the domain D of M are necessarily definable in P. We consider the significance of introducing such a restriction in the next section.

We note that the explicit introduction of the concept of a "terminating routine" highlights an ambiguity[96] that is implicit in classical definitions. Thus, Definition 5(*iv*) can be understood to mean either of[97]:

**Definition 5**(*iv*)(*a*): If [*R*] is a formula of P, *a sequence s of* M *satisfies* [(A$x_i$)*R*] *individually*, under the interpretation M, if, and only if, given any sequence *s'* of M, which differs from *s* in at most the *i*'th component, there is an *individually terminating routine* to determine that *s'* satisfies [*R*] classically.

**Definition 5**(*iv*)(*b*): If [*R*] is a formula of P, *a sequence s of* M *satisfies* [(A$x_i$)*R*] *uniformly*, under the interpretation M, if, and only if, there is a *uniformly terminating routine* to determine that, given any sequence *s'* of M, which differs from *s* in at most the *i*'th component, *s'* satisfies [*R*] classically.

We note, here, that:

**Meta-lemma 5**: 5(*iv*)(*b*) ==> 5(*iv*)(*a*) ==> 5(*iv*)[98]

*Proof*: Clearly, if there is a uniformly terminating routine to determine that, given any sequence *s'* of M, which differs from *s* in at most the *i*'th component, *s'* satisfies [*R*] classically, then, trivially, given any sequence *s'* of M, which differs from *s* in at most the *i*'th component, the routine determines that *s'* satisfies [*R*] classically. Hence, if

---

[96] We note that classical set theory is marked by a similar ambiguity; it also does not explicitly distinguish between a mapping that is determined by an individually terminating routine, and a mapping that is determined by a uniformly terminating routine. Since we may have an individually terminating routine to determine that two functions are equivalent in their instantiations, but no uniformly terminating routine that also establishes such equivalence, it follows, from Corollary 1.2 to Meta-theorem 1, that the classical definition of a function as a mapping may lead to anomalies through such ambiguity.

[97] In other words, a routine is a terminating routine if, and only if, it is either an individually terminating routine, or a uniformly terminating routine.

[98] We use this shorthand notation to mean that if a sequence is satisfied uniformly under M, then it is satisfied individually under M; and that the latter implies it is satisfied classically under M.

5(*iv*)(*b*) holds, then 5(*iv*)(*a*) holds trivially. Since 5(*iv*) implies that either 5(*iv*)(*a*) or 5(*iv*)(*b*) holds, the meta-lemma follows.¶

**Definition 5**(*v*): A formula [*R*] of P is *classically true*, for the interpretation M, if, and only if, there is a terminating routine to determine that every sequence of M satisfies [*R*] classically.

As in 5(*iv*) above, 5(*v*) can be taken to mean either of:

> **Definition 5**(*v*)(*a*): A formula [*R*] of P is *individually true*, for the interpretation M, if, and only if, given any sequence *s* of M, there is an individually terminating routine to determine that *s* satisfies [*R*] classically.

> **Definition 5**(*v*)(*b*): A formula [*R*] of P is *uniformly true*, for the interpretation M, if, and only if, there is a uniformly terminating routine to determine that, given any sequence *s* of M, *s* satisfies [*R*] classically.

We again have:

> **Meta-lemma 6**: 5(*v*)(*b*) ==> 5(*v*)(*a*) ==> 5(*v*)

> *Proof*: As in the previous meta-lemma, if 5(*v*)(*b*) holds, then 5(*v*)(*a*) holds trivially. Since 5(*v*) implies that either 5(*v*)(*a*) or 5(*v*)(*b*) holds, the meta-lemma follows.¶

**Definition 5**(*vi*): A formula [*R*] of P is *classically false*, for the interpretation M, if, and only if, there is a terminating routine to determine that no sequence of M satisfies [*R*] classically.

Here, too, (*vi*) can be taken to mean either of:

**Definition 5**(*vi*)(*a*): A formula [*R*] of P is *individually false*, for the interpretation M, if, and only if, given any sequence *s* of M, there is an individually terminating routine to determine that *s* does not satisfy [*R*] classically.

**Definition 5**(*vi*)(*b*): A formula [*R*] of P is *uniformly false*, for the interpretation M, if, and only if, there is a uniformly terminating routine to determine that, given any sequence *s* of M, *s* does not satisfy [*R*] classically.

We also have:

**Meta-lemma 7**: 5(*vi*)(*b*) ==> 5(*vi*)(*a*) ==> 5(*vi*)

*Proof*: As in the previous meta-lemma, if 5(*vi*)(*b*) holds, then 5(*vi*)(*a*) holds trivially. Since 5(*vi*) implies that either 5(*vi*)(*a*) or 5(*vi*)(*b*) holds, the meta-lemma follows.¶

## 5.2  A constructive expression of Church's Thesis

Now we note, firstly, that definitions 5(*iv*)(*a*), 5(*v*)(*a*) and 5(*vi*)(*a*) are meta-assertions that appeal to the decidability of an infinity of individual assertions in M, between elements of D that are not necessarily representable in P as mathematical objects. In the absence of a terminating routine for deciding whether an assertion holds individually in M or not, the definitions are, essentially, non-constructive.

We note, next, that definitions 5(*iv*)(*b*), 5(*v*)(*b*) and 5(*vi*)(*b*) are also meta-assertions that appeal to the decidability of infinitely compound assertions in M, between elements of D that are not necessarily representable in P as mathematical objects. In the absence of a terminating routine for deciding whether an infinitely compound assertion is satisfied uniformly in M or not, these definitions, too, are, essentially, non-constructive.

We note also that, if the domain of M is representable in P, then definitions 5(*i*) to 5(*iii*) can be treated as finitely decidable "truths of fact". These definitions are, then, constructive, and intuitionistically unobjectionable. The constructiveness of definitions 5(*iv*) to 5(*vi*), however, cannot be assumed even in such a case.

The questions arise:

(*) When may we constructively assume that, given any sequence *s* of an interpretation M of P, there is an individually terminating routine to determine that *s* satisfies a given P-formula in M?

(**) When may we constructively assume that there is a uniformly terminating routine such that, given any sequence *s* of an interpretation M, *s* satisfies a given P-formula in M?

We note that, if the domain D of M can be assumed representable in P, then (*) can be answered constructively. More precisely, we reformulate the classical Church Thesis[99]:

**Individual Church Thesis**: If, for a given relation $R(x_1, ..., x_n)$, and any sequence $<a_1, ..., a_n>$, in some interpretation M of P, there is an individually terminating routine such that it will determine whether $R(a_1, ..., a_n)$ holds in M or not, then every element of the domain D of M is the interpretation of some term of P, and there is some P-formula $[R'(x_1, ..., x_n)]$ such that:

$R(a_1, ..., a_n)$ holds in M if, and only if, $[R'(a_1, ..., a_n)]$ is P-provable.

In other words, the Individual Church Thesis postulates that, if a relation *R* is effectively decidable individually (which may be non-algorithmically) in an interpretation M of some formal system P, then *R* is expressible in P, and its domain necessarily consists of only mathematical objects, even if the predicate letter *R* is not, itself, a mathematical object.

[99] The relation between the classical Church Thesis and the reformulated Theses is brought out by Corollary 14.2.

However, (\*\*) can be answered constructively for any interpretation M of P, if we postulate:

**Uniform Church Thesis**: If, in some interpretation M of P, there is a uniformly terminating routine such that, for a given relation $R(x_1, ..., x_n)$, and any sequence $<a_1, ..., a_n>$, it will determine whether $R(a_1, ..., a_n)$ holds in M or not, then $R(x_1, ..., x_n)$ is the interpretation in M of a P-formula $[R(x_1, ..., x_n)]$, and:

$R(a_1, ..., a_n)$ holds in M if, and only if, $[R(a_1, ..., a_n)]$ is P-provable.

Thus, the Uniform Church Thesis postulates that, if a relation *R* is effectively decidable uniformly (necessarily algorithmically) in an interpretation M of a formal system P, then, firstly, *R* is expressible in P, and, secondly, the predicate letter *R*, and all the elements in the domain of the relation *R*, are necessarily mathematical objects.

It follows from definition 5(*v*)(*b*) that:

**Meta-lemma 8**: The Uniform Church Thesis implies that a formula $[R]$ is P-provable if, and only if, $[R]$ is uniformly true in some interpretation M of P.

**Proof**: By definition 5(*v*)(*b*), if $[R]$ is uniformly true in some interpretation M of P, then there is some uniformly terminating routine such that, for a given relation $R(x_1, ..., x_n)$, and any sequence $<a_1, ..., a_n>$, it will determine whether $R(a_1, ..., a_n)$ holds in M or not. By the Uniform Church Thesis, it follows that:

$R(a_1, ..., a_n)$ holds in M if, and only if, $[R(a_1, ..., a_n)]$ is P-provable.

Hence, the domain of M is representable in P, and therefore denumerable. It follows that the domains of any two interpretations of P are isomorphic ([Me64], p90), and, by definition 5(*v*)(*b*), [*R*] is uniformly true in any given interpretation M of P. By Gödel's Completeness Theorem ([Me64], p68, Corollary 2.15(a)), [*R*] is P-provable.¶

**Corollary 8**.1: The Uniform Church Thesis implies that, if there is a uniformly terminating routine such that the number-theoretic relation *R* is satisfied by any sequence *s* in some interpretation M of P, then the predicate letter "*R*" is a formal mathematical object in P.

We also have the further meta-lemma:

**Meta-lemma 9**: The Uniform Church Thesis implies that, if a P-formula [*R*] is uniformly true in some interpretation M of P, then [*R*] is uniformly true in every model of P.

**Corollary 9.1**: The Uniform Church Thesis implies that if a formula [*R*] is not P-provable, but [*R*] is classically true under the standard interpretation, then [*R*] is individually true, but not uniformly true, in the standard model of P.

**Corollary 9.2**: The Uniform Church Thesis implies that Gödel's undecidable sentence ([Go31a], p26) is individually true, but not uniformly true in the standard model of P.

## 6. Constructive definitions of classical concepts

As suggested by the constructive expression of Church's Thesis, a more stringent condition of constructivity would be to postulate that assertions in an interpretation M of P are verifiable individually only over elements in the domain D of M that are representable in P as mathematical objects. We thus have the following constructive definitions[100], which correspond to various classical concepts that involve individually terminating routines.

**Definition 6**(*i*): A relation $R(x_1, ..., x_n)$ of M is *effectively expressible* in P if, and only if, there is a P-formula $[R'(x_1, ..., x_n)]$, with $n$ free variables, such that, for any sequence of numerals $<[b], [a_1], ..., [a_n]>$ of P, there is an individually terminating routine for determining that:

    (*i*)   if $R(a_1, ..., a_n)$ holds in M, then $[R'(a_1, ..., a_n)]$ is P-provable;

    (*ii*)  if $R(a_1, ..., a_n)$ does not hold in M, then $[\sim R'(a_1, ..., a_n)]$ is P-provable,

    where $<b, a_1, ..., a_n>$ is the interpretation of $<[b], [a_1], ..., [a_n]>$ in M.

We note that, for any given sequence of numerals $<[b], [a_1], ..., [a_n]>$ of P, in the absence of an individually terminating routine for determining (*i*) and (*ii*), the concept of "effective expressibility" is, essentially, non-constructive.

**Definition 6**(*ii*): A function $F(x_1, ..., x_n)$ of M is *effectively representable* in P if, and only if, there is a P-formula $[R(y, x_1, ..., x_n)]$, with the free variables $[(y, x_1, ..., x_n)]$, such that, for any sequence of numerals $<[b], [a_1], ..., [a_n]>$ of P:

    (*i*)   there is an individually terminating routine for determining that if ($F(a_1, ..., a_n)$ $= b$) holds in M, then $[R(b, a_1, ..., a_n)]$ is P-provable;

    (*ii*)  $[(E!y)R(y, a_1, ..., a_n)]$ is P-provable,

---

[100] The following definitions are based on extending the classical definitions of the "expressibility", "representability" and "strong representability" of number-theoretic relations and functions in a first-order theory, as expressed by Mendelson ([Me64], p118).

where $<b, a_1, ..., a_n>$ is the interpretation of $<[b], [a_1], ..., [a_n]>$ in M.

We again note that, for any given sequence of numerals $<[b], [a_1], ..., [a_n]>$ of P, in the absence of an individually terminating routine for determining (*i*), the concept of "effective representability" is, essentially, non-constructive.

**Definition 6**(*iii*): A function $F(x_1, ..., x_n)$ of M is *strongly representable effectively* in P if, and only if, there is a P-formula $[R(y, x_1, ..., x_n)]$, with the free variables $[(y, x_1, ..., x_n)]$, such that:

(*i*)   for any sequence of numerals $<[b], [a_1], ..., [a_n]>$ of P, there is an individually terminating routine for determining that if $(F(a_1, ..., a_n) = b)$ holds in M, then $[R(b, a_1, ..., a_n)]$ is P-provable;

(*ii*)  $[(E!y)R(y, x_1, ..., x_n)]$ is P-provable,

where $<b, a_1, ..., a_n>$ is the interpretation of $<[b], [a_1], ..., [a_n]>$ in M.

We note that, here also, for any given sequence of numerals $<[b], [a_1], ..., [a_n]>$ of P, in the absence of an individually terminating routine for determining (*i*), the concept of "effectively strong representability" is, essentially, non-constructive.

**Definition 6**(*iv*): A relation $R(x_1, ..., x_n)$ in M is *effectively definable* in P if, and only if, $[R(x_1, ..., x_n)]$ is a P-formula such that, for any sequence of numerals $<[a_1], ..., [a_n]>$ of P, there is an individually terminating routine for determining that:

(*i*)   if $R(a_1, ..., a_n)$ holds in M, then $[R(a_1, ..., a_n)]$ is P-provable;

(*ii*)  if $R(a_1, ..., a_n)$ does not hold in M, then $[\sim R(a_1, ..., a_n)]$ is P-provable,

where $<a_1, ..., a_n>$ is the interpretation of $<[a_1], ..., [a_n]>$ in M.

**Meta-lemma 10**: If a relation $R$ of M is effectively definable in P, then it is the interpretation of one of its representations in P. (Although the converse holds if M is the standard model, it need not be true for every model.)

We note that, for any given sequence of numerals $<[b], [a_1], ..., [a_n]>$ of P, in the absence of an individually terminating routine for determining (*i*) and (*ii*), the concept of "effective definability" is also, essentially, non-constructive.

**Definition 6**(*v*)**:** A function $F(x_1, ..., x_n)$ is *effectively computable individually* in an interpretation M of P if, and only if, it is effectively representable in P.

**Definition 6**(*vi*)**:** A function $F(x_1, ..., x_n)$ is *effectively computable uniformly* in an interpretation M of P, if, and only if, it is both effectively representable and effectively definable in P.

**Meta-lemma 11**: If a function $F(x_1, ..., x_n)$ is effectively computable uniformly in an interpretation M of P, then it is effectively computable individually in M; however, the converse is not true.

*Proof*: The first part is trivially true, by definition. The converse follows from Meta-theorem 1.¶

**Definition 6**(*vii*)**:** A function $F(x_1, ..., x_n)$ is an *omega-function* in an interpretation M of P if, and only if, there is a P-formula $[R(y, x_1, ..., x_n)]$, with the free variables $[(y, x_1, ..., x_n)]$, such that:

(*i*)   for any sequence of numerals $<[b], [a_1], ..., [a_n]>$ of P, there is an individually terminating routine for determining that, if $(F(a_1, ..., a_n) = b)$ holds in M, then $[R(b, a_1, ..., a_n)]$ is P-provable, where $<b, a_1, ..., a_n>$ is the interpretation of $<[b], [a_1], ..., [a_n]>$ in M;

(*ii*)  for any sequence of numerals $<[b], [a_1], ..., [a_n]>$ of P, there is an individually terminating routine for determining that $[(E!y)(R(b, a_1, ..., a_n))]$ is P-provable;

(*iii*) $[(E!y)(R(y, x_1, ..., x_n))]$ is not-P-provable,

where $<b, a_1, ..., a_n>$ is the interpretation of $<[b], [a_1], ..., [a_n]>$ in M.

We thus have:

**Meta-lemma 12**: An omega-function $F(x_1, ..., x_n)$ in an interpretation M of P is effectively computable individually, but not effectively computable uniformly.

We can now define:

**Definition 6**(*viii*): A number-theoretic function $F(x_1, ..., x_n)$ in the standard interpretation M of P is *uncomputable* if, and only if, it is an omega-function.

## 7. Self-terminating, converging and oscillating Turing machines

We note that classical Turing machines ([Me64], p229) are constructive, in the sense of

§II-6, to the extent that:

**Meta-lemma 13**: If a total function $F(x_1, ..., x_n)$ is classically Turing-computable[101]

in the standard interpretation M of P, then it is effectively computable individually in

M.[102]

---

[101] We call a function $F$ "classically Turing-computable" if, and only if, there is a Turing-computable algorithm $U$ that computes $F$ as defined by Mendelson ([Me64], p231). See §II-7(1)(*l*) below.

[102] We note that the classical Turing Thesis is essentially the assertion that if a function $F(x_1, ..., x_n)$ is effectively computable individually in the standard interpretation M of P, then it is classically Turing-computable in M.

*Proof*: Since the domain of M is representable in P, a number-theoretic function is representable (cf. [Me64], p118) in P if, and only if, it is effectively representable in P (cf. §II-6(*ii*)). Now, if a total number-theoretic function $F(x_1, ..., x_n)$ is classically Turing-computable in the standard interpretation M of P, then it is recursive ([Me64], p233, Corollary 5.13). Since every recursive function is representable in P ([Me64], p131, Proposition 3.23), it follows by Def 5 (§II-6) that $F(x_1, ..., x_n)$ is effectively computable individually in M.¶

We next define a neo-classical Turing machine NT as a natural, and essentially constructive, extension of a classical Turing machine T. We introduce the concept of a neo-classical, self-terminating, Turing routine, and note some significant consequences of the difference between such routines and classical, program-terminating[103], Turing routines. We begin by recalling Mendelson's description ([Me64], p229) of the operations of T:

### 7.1 Classical Turing machines

(*a*)   There is a two-way, potentially infinite[104], tape divided up into squares. There is a finite set of *tape symbols* $S_0, S_1, ..., S_n$ called the *alphabet* of the machine; at every moment, each square of the tape is occupied by at most one symbol. The machine has a finite set of *internal states* $\{q_0, q_1, ..., q_m\}$. At any given moment, the machine is in exactly one of these states. Finally, there is a reading head which, at any given time, stands over some square of the tape. The machine does not act continuously, but only at discrete moments of time. If, at any moment $t$, the reading head is scanning a square containing a symbol[105] $S_i$, and the machine is in the internal state $q_j$, then the action of the machine is determined, and it will do one of four things:

> (*i*)   it may erase the symbol $S_i$ and print a new symbol $S_k$;

> (*ii*)  it may move left one square;

> (*iii*) it may move right one square;

> (*iv*)  it may stop. In cases (*i*)-(*iii*), the machine goes into a new internal state $q_r$, and is ready to act again at time $t$+1.

---

[103] Note that only what we term as a "program-terminating routine" is classically defined as a "halting routine"; all others are termed as "non-halting routines".

[104] "The tape is said to be potentially infinite in the sense that, although at any moment it is finite in length, additional squares can be added to the right- and left-hand ends of the tape." ([Me64], p229).

(*b*)   The actions (*a*)(*i*) to (*a*)(*iii*) can be represented by *quadruples*:

> (*i*)   $q_j S_i S_k q_r$;

> (*ii*)  $q_j S_i L q_r$;

> (*iii*) $q_j S_i R q_r$,

where the first two symbols stand for the present internal state and scanned symbol, the third symbol represents the action of the machine, and the fourth symbol gives the internal state of the machine after the action has been performed.

(*c*)   If a tape is put into a Turing machine and the reading head is placed on a certain square, and if the machine is started off in one of its internal states, then the machine begins to operate on the tape: printing and erasing symbols and moving from one square to an adjacent one. If the machine ever stops, the resulting tape is said to be the *output* of the machine applied top the given tape.

---

[105] "We shall assume that the symbol $S_0$ represents a blank, so that the reading head may always be assumed to be scanning a symbol." ([Me64], p229)

(*d*)  We can associate with any Turing machine T the following algorithm *B* in the alphabet A of T. Take any word *P* of the alphabet A and print it from left to right in the squares of an empty tape. Place this tape in the machine with the reading head scanning the left-most square. Start the machine in the internal state $q_0$. If the machine ever stops, the word of A appearing on the tape is the value of the algorithm *B*. *B* is called a *Turing algorithm*.

(*e*)  A *classical Turing machine* T is then defined precisely as a finite set of quadruples of the three kinds (*b*)(*i*) to (*b*)(*iii*), where no two quadruples have the same first two symbols.

(*f*)  An *instantaneous tape description* of a Turing machine T is a word such that:

  (*i*)   all symbols in the word but one are tape symbols $S_m$;

  (*ii*)  the only symbol which is not a tape symbol is an internal state $q_s$;

  (*iii*) $q_s$ is not the last symbol of the word.[106]

(*g*)  We say that T *moves* one instantaneous tape description *alpha* into another one *beta*[107] if, and only if:

  (*i*)   either, *alpha* is of the form $Pq_jS_iQ$, *beta* is of the form $Pq_rS_kQ$, and $q_jS_iS_kq_r$ is one of the quadruples of T;

---

[106] "An instantaneous tape description describes the condition of the machine and the tape at a given moment. When read from left to right, the tape symbols in the description represent the symbols on the tape at the moment. The internal state $q_s$ in the description  is the internal state of the machine at the moment, and the tape symbol occurring immediately to the right of $q_s$ in the tape description represents the symbol being scanned by the machine at the moment." ([Me64], p230, footnote 1).

[107] Abbreviated by "T: *alpha -> beta*".

(*ii*) or, *alpha* is of the form $PS_sq_jS_iQ$, *beta* is $Pq_rS_sS_iQ$, and $q_jS_iLq_r$ is one of the quadruples of T;

(*iii*) or, *alpha* is of the form $q_jS_iQ$, *beta* is $q_rS_0S_iQ$, and $q_jS_iLq_r$ is one of the quadruples of T;

(*iv*) or, *alpha* is of the form $Pq_jS_iS_kQ$, *beta* is $PS_iq_rS_kQ$, and $q_jS_iRq_r$ is one of the quadruples of T;

(*v*) or, *alpha* is of the form $Pq_jS_i$, *beta* is $PS_iq_rS_o$, and $q_jS_iRq_r$ is one of the quadruples of T.[108]

(*h*) We say that T *stops*[109] *classically* at an instantaneous tape description *alpha* if, and only if, there is no instantaneous tape description *beta* into which T can move *alpha*.

(*j*) A *classical computation* of a Turing machine T is a finite sequence of instantaneous tape descriptions $alpha_0$, $alpha_1$, ..., $alpha_m$ ($m >= 0$) such that:

(*i*) the internal state occuring in $alpha_0$ is $q_0$;

(*ii*) for $0 =< i < m$, $alpha_{(i+1)}$ follows $alpha_i$;

(*iii*) there is no instantaneous tape description $alpha_{(m+1)}$ into which T can move $alpha_m$;

(*iv*) and T stops at $alpha_m$.

---

[108] "Observe that, according to our intuitive picture, 'T moves *alpha* into *beta*' implies that if the condition at time $t$ of the Turing machine and tape is described by *alpha*, then the condition at time $t+1$ is described by *beta*. Notice that, according to (*g*)(*iii*), whenever the machine reaches the left-hand end of the tape and is ordered to move left, a blank square is attached to the tape on the left; similarly, by (*g*)(*v*), a blank square is added on the right when the machine reaches the right-hand end of the tape and has to move right." ([Me64], p230, footnote 2).

[109] We note that the terms "stops" and "halts" are synonymous.

(*k*)   The algorithm $B_{T,C}$ in any alphabet C containing the alphabet A of T is defined as follows:

For any words $P$, $Q$ in C, $B_{T,C}(P) = Q$ if, and only if, there is a classical computation of T which begins with the  instantaneous tape description $q_0P$ and ends with an instantaneous tape description of the form $R_1q_jR_2$, where $Q = R_1R_2$.

(*l*)   An algorithm $U$ in an alphabet D is called *classically Turing-computable* if, and only if, there is a Turing machine T with alphabet A and an alphabet C containing A+D[110] such that $B_{T,C}$ and $U$ are fully equivalent relative to D.

(*m*)   Given a partial number-theoretic function $F(x_1, ..., x_n)$, we say that a Turing machine T (whose alphabet A includes {1, *}) computes $F$ if, and only if, for any natural numbers $k_1$, $k_2$, ... $k_n$, and any word $Q$, $B_{T,A}(k_1* k_2* ... * k_n) = Q$ if, and only if, $Q$ is $R_1F(x_1, ..., x_n)R_2$, where both $R_1$ and $R_2$ are certain, possibly empty, words consisting only of $S_0$'s.[111]

## 7.2  Neo-classical Turing machines

We now define a broader class of neo-classical Turing machines NT as follows:

(*a*)   A neo-classical Turing machine NT is a classical Turing machine T that also stops if an instantaneous tape description repeats itself.[112]

---

[110] By "A+D" we mean the combined alphabet of A and D.

[111] We note that $[k_1* k_2*...k_n]$ represents, in A, the $k$-tuple $k_1$, $k_2$, ..., $k_n$, ([Me64], p212), and that $S_0$ is interpreted as a blank ([Me64], p231).

[112] It is convenient to visualise a neo-classical Turing machine as a two-dimensional virtual-teleprinter, which maintains a copy of every instantaneous tape description in a random-access memory during a computation.

(*b*) A *self-terminating computation* of an NT machine is a finite sequence of instantaneous tape descriptions $alpha_0$, $alpha_1$, ..., $alpha_m$ ($m >= 0$) such that:

    (*i*)   the internal state occuring in $alpha_0$ is $q_0$;

    (*ii*)  for $0 =< i < m$, $alpha_{(i+1)}$ follows $alpha_i$;

    (*iii*) $alpha_m$ is a repetition of $alpha_i$ for some $0 =< i < m$;

    (*iv*) NT stops at $alpha_m$.[113]

(*c*) We call classical Turing-computations as *program-terminating* NT-*routines*, and neo-classical Turing-computations as *terminating* NT-*routines*.

(*d*) The *instantaneous tape output* of an NT machine is the word obtained by deleting the symbol for the internal state of the machine from the instantaneous tape description.

(*e*) The number of tape symbols in an instantaneous tape output of an NT machine is called the *output length* of the instantaneous tape description.

(*f*) For all $0 =< i =< m$, the finite integer $i$ is defined as the *configuration number* of an *instantaneous tape output* corresponding to the instantaneous tape description $alpha_i$ of an NT machine.

---

[113] We assume that a neo-classical Turing machine contains some effective terminating routine for comparing instantaneous tape descriptions after executing an instruction defined by a quadruple, and for halting when an instantaneous tape description repeats itself. We also assume that the alphabet of such a machine includes a special symbol for self-termination that is returned in such a case.

We note that such a symbol would not occur in any quadruple; strictly speaking, it would be treated more appropriately as a meta-symbol that determines a meta-action of an NT machine. We may, indeed, reasonably argue that an NT machine is essentially, if not precisely, the appropriate parallel of a meta-proof (such as, for instance, Gödel's argument that determines a formula $[F(n)]$ as P-provable for every numeral $[n]$, even though $[F(x)]$ is not P-provable).

(*g*) The *terminal tape output* of an NT machine is the word obtained by deleting the symbol for the internal state of the machine from the final instantaneous tape description of a terminating NT-routine.

(*h*) An NT-routine is *non-terminating* if, and only if, it is not a terminating NT-routine.

(*j*) A non-terminating NT-routine is a *converging* NT-*routine* if, and only if, there is a positive integer $n_0$ such that, given any natural number $n > n_0$, there is a minimum configuration number $n_{min}$ such that, for any given configuration number $m > n_{min}$, the instantaneous tape description is of the form $Pq_jS_iQ$ where $P$ is a word in the alphabet A of NT that contains more than $n$ tape symbols.

(*k*) We call a non-terminating NT-routine an *oscillating* NT-*routine* if, and only if, it is not a *converging* NT-*routine*.

### 7.3 Total partial recursive functions

The significance, of defining the truth of a formula of P under interpretation explicitly in terms of terminating routines, of expressing Church's Thesis constructively, and of defining self-terminating computations of an NT machine, is expressed by the following meta-lemmas.

**Meta-lemma 14**: If we assume a Uniform Church Thesis, then every partial recursive[114] number-theoretic function $F(x_1, ..., x_n)$ has a unique constructive extension as a total function.

---

[114] Classically ([Me64], p120-121, p214), a partial function $F$ of $n$ arguments is called partial recursive if, and only if, $F$ can be obtained from the initial functions (zero function), projection functions, and successor function (of classical recursive function theory) by means of substitution, recursion and the classical, unrestricted, $\mu$-operator. $F$ is said to come from $G$ by means of the unrestricted $\mu$-operator, where

*Proof*: We assume that $F$ is obtained from the recursive function $G$ by means of the unrestricted $\mu$-operator, so that $F(x_1, ..., x_n) = \mu y(G(x_1, ..., x_n, y) = 0)$.

If $[H(x_1, ..., x_n, y)]$ expresses $\sim(G(x_1, ..., x_n, y) = 0)$ in P, we consider the P-provability, and truth in the standard interpretation M of P, of the formula $[H(a_1, ..., a_n, y)]$ for a given sequence of numerals $<[a_1], ..., [a_n]>$ of P.

(We note that, if we define the truth of a formula of P, under an interpretation M, explicitly in terms of terminating routines that are appropriate to M[115], a formula such as $[(Ax)F(x)]$ does not interpret as a compound number-theoretic assertion about the range of values of M for which $[F(x)]$ is satisfied in M, but as a meta-assertion that there is a terminating routine, that is appropriate to M, which can determine that $F(a)$ holds in M for any given $a$ of M.)

We now consider the argument:

(*a*) Let $Q_1$ be the meta-assertion that $[H(a_1, ..., a_n, y)]$ is not classically true in M. Hence there is no terminating routine in M such that, for any given $y$ in M, $y$ satisfies $[H(a_1, ..., a_n, y)]$ classically. It follows that there is no uniformly terminating routine in M such that, for any given $y$ in M, $y$ satisfies $[H(a_1, ..., a_n, y)]$ classically.

---

$G(x_1, ..., x_n, y)$ is recursive, if, and only if, $F(x_1, ..., x_n) = \mu y(G(x_1, ..., x_n, y) = 0)$, where $\mu y(G(x_1, ..., x_n, y) = 0)$ is the least number $k$ (if such exists) such that, if $0 = <i = <k$, $G(x_1, ..., x_n, i)$ exists and is not 0, and $G(x_1, ..., x_n, k) = 0$. We note that, classically, $F$ may not be defined for certain $n$-tuples; in particular, for those $n$-tuples $(x_1, ..., x_n)$ for which there is no $y$ such that $G(x_1, ..., x_n, y) = 0$. (We note that the classical $\mu$-operator ([Me64], p121), and the *e*-operator ([Go31a], p16), defined in Meta-theorem 1(*iv*), are identical.)

[115] We note that a terminating routine that is appropriate to M need not be finite, or even denumerable, unless it is specified as constructive in the sense of §II-6.

Since $G(a_1, ..., a_n, y)$ is recursive, it follows that there is some finite $k$ such that any neo-classical Turing machine $NT_1(y)$ that computes $G(a_1, ..., a_n, y)$ will halt and return the value 0 for $y = k$.

(b)  Next, let $Q_2$ be the meta-assertion that $[H(a_1, ..., a_n, y)]$ is classically true in the standard interpretation M of P, but that there is no uniformly terminating routine in M such that, for any given $y$ in M, $y$ satisfies $[H(a_1, ..., a_n, y)]$ classically.

Since $G(a_1, ..., a_n, y)$ is recursive, it follows that there is some finite $k$ such that the neo-classical Turing machine $NT_1(y)$ will halt, and return the symbol for self-termination for $y = k$.

(c)  Finally, let $Q_3$ be the meta-assertion that $[H(a_1, ..., a_n, y)]$ is classically true in the standard interpretation M of P, and that there is a uniformly terminating routine in M such that, for any given $y$ in M, $y$ satisfies $[H(a_1, ..., a_n, y)]$ classically. We then have that that $[H(a_1, ..., a_n, y)]$ is uniformly true in the standard interpretation M of P

If we assume a Uniform Church Thesis, then, by Meta-lemma 8, it follows that $[H(a_1, ..., a_n, y)]$ is P-provable. Let $h$ be the Gödel-number of $[H(a_1, ..., a_n, y)]$. We consider, then, Gödel's recursive number-theoretic relation $xBy$, which holds in M if, and only if, $x$ is the Gödel-number of a proof sequence in P for the P-formula whose Gödel-number is $y$. It follows that there is some finite $k$ such that any neo-classical Turing machine $NT_2(y)$, which computes the characteristic function[116] of $xBh$, will halt and return the value 0 for $x = k$[117].

---

[116] "If $R(x_1, ..., x_n)$ is a relation, then the characteristic function $C_n(x_1, ..., x_n)$ is defined as follows:

$C_n(x_1, ..., x_n) = 0$ if $R(x_1, ..., x_n)$ is true, and $C_n(x_1, ..., x_n) = 1$ if $R(x_1, ..., x_n)$ is false." ([Me64], p119).

Since $Q_1$, $Q_2$, and $Q_3$ are mutually exclusive and exhaustive[118], it follows that, when run simultaneously[119] over the sequence 1, 2, 3, ... of values for $y$, one of the parallel duo $\{NT_1(y) // NT_2(y)\}$ will always halt for some finite value of $y$. If $NT_1(y)$ halts at $y = k$, and returns the value 0, we define $F'(a_1, ..., a_n) = F(a_1, ..., a_n)$. If $NT_1(y)$ halts and returns the symbol for self-termination, or if $NT_2(y)$ halts, we define $F'(a_1, ..., a_n) = 0$.

Hence, under the given hypothesis, there is always a unique constructive extension of every partial recursive function $F(x_1, ..., x_n)$ as a total function $F'(x_1, ..., x_n)$.¶

---

[117] We assume that such a machine can be effectively meta-programmed to proceed to the next instantaneous tape description whenever it encounters a loop.

[118] They correspond to the instances where a classical Turing machine that computes the recursive function $G(a_1, ..., a_n, y)$ will halt for some $y$, loop for some $y$, or not halt for any $y$, respectively.

[119] This concept is essentially that of parallel computing, where the action of one machine can influence the action of another unpredictably, without human intervention. Since classical Turing machines are necessarily sequential, such a procedure cannot be defined as a classical Turing machine. In his article "Uncomputability in the work of Alan Turing and Roger Penrose" (a talk given for Interface 5, Hamburg, 6 October 2000, and available at <http://www.turing.org.uk/philosophy/lecture1.html>), Andrew Hodges remarks that the possibility of parallel machines being essentially different from his Logical Computing Machines does not (arguably) appear to have been considered by Turing:

"... Another source may lie in Turing's definition of an 'oracle-machine' which is a Turing machine allowed at certain points to 'consult the oracle'. Such a machine is not purely mechanical: it is like the 'choice-machine' defined in (Turing 1936-7) which at certain points allows for human choices to be made. Turing used the word 'machine' for entities which are only partially mechanical in operation, reserving the term 'automatic machine' for those which are purely mechanical. Copeland appears to imagine that when Turing describes the oracle-machine definition as giving a 'new type of machine', he is defining a new type of automatic machine. On the contrary, Turing is defining something only partially mechanical.

To take this point further, it is worth noting that the expression 'purely mechanical process' enters into Turing's definitive statement of the Church-Turing thesis, which comes as an opening section to (Turing 1939), and that Turing goes on: 'understanding by a purely mechanical process one which could be carried out by a machine'. In the subsequent discussion the word 'machine' is used to mean 'Turing machine'. There is no evidence that Turing had any concept of a purely mechanical 'machine' of any kind other than encapsulated by the Turing machine definition."

It now follows from Meta-lemma 14:

**Corollary 14.1**: The classical Halting problem is effectively solvable if we assume a Uniform Church Thesis.

## 7.4  The Uniform Church Thesis and the classical Church-Turing Theses

We now consider the significance of constructively defining terminating routines for the two classical theses:

**Classical Church's Thesis**: A number-theoretic function is effectively computable (partially) if, and only if, it is (partially) recursive ([Me64], p147, p227).

**Classical Turing Thesis**: Every effectively computable function is classically Turing computable ([Me64], p237).

We note that, classically, the two theses are accepted as equivalent since the Turing machine approach to effective computability is considered equivalent[120] ([Me64], p237) to that by means of normal algorithms ([Me64], p209), or by recursive functions.

Now, we can express the Turing Thesis alternatively in terms of terminating routines:

**Alternate Turing Thesis**: A function is classically Turing computable if, and only if, it is computable[121] either by an individually terminating routine, or by a uniformly terminating routine.

---

[120] However, such equivalence is based on the argument that there is a uniformly effective method (algorithm) for computing any (partial) recursive function. The classical proof of this uses induction over (partial) recursive functions, thus assuming that they are mathematical objects that have properties that are within the scope of mathematical induction. By Meta-lemma 1, such an assumption is invalid.

[121] In the sense of §II-4(c).

We also note that the classical Turing Thesis - that every effectively computable function is classically Turing computable ([Me64], p237) - implies the:

**Secondary Turing Thesis**: If a number-theoretic function $F(x_1, ..., x_n)$ is effectively computable individually in the standard interpretation M of P, then it is classically Turing computable in M.

However, if every partial recursive function can be constructively extended as a total function, which is effectively computable individually in the standard interpretation M of P by the argument in Meta-lemma 14, then it follows that the classical Turing Thesis is false, since:

**Corollary 14.2**: If we assume a Uniform Church Thesis, then not every effectively computable function is classically Turing computable.[122]

**Proof:** Meta-lemma 14 gives an effective method, corresponding to a constructive Turing "oracle"[123], for determining whether a classical Turing machine will halt or

---

[122] We thus have the distinct possibility that adopting the classical Church-Turing Thesis may limit the ability of our mathematical language to adequately express mathematical aspects - pertaining to the subjects of our intuitive experience - that are not necessarily Turing-computable.

[123] We note that this constructivity follows from reformulating Church's Thesis constructively. In contrast, based on the non-constructivity implicit in his own thesis regarding effective computability, Turing appears to have believed that his concept of a mathematical "oracle" must remain essentially non-constructive. As described by Hodges (op cit):

"However the driving force lay in the question: what is the consequence of supplementing a formal system with uncomputable deductive steps? In pursuit of this question, Turing introduced the definition of an 'oracle' which can supply on demand the answer to the halting problem for every Turing machine. Turing gave his subject-matter an interpretation which described the mathematician's 'intuition' in theorem-proving, and Newman (1955) effectively identified the uncomputable 'oracle' with intuition. This was perhaps going too far since the 'oracle' is capable of far more than any human being; nevertheless Newman had a unique status as Turing's collaborator at this period and must reflect the tenor of Turing's discussions. In any case, Turing makes it clear that the 'intuition' being discussed is related to the human act of seeing the truth of a formally unprovable Gödel statement. To summarise, it is notable that Turing's 1938 work focussed on the same issue as Penrose now raises: the interpretation of uncomputable deductions.

not. Hence, by the classical Halting argument, it does not define a classical Turing machine.¶

Since the Individual Church Thesis implies the classical Church's Thesis, it further follows that:

**Corollary 14.3**: If we assume a Uniform Church Thesis, then not every (partially) recursive function is classically Turing-computable.[124]

## 7.5  Converging NT-routines and Cauchy sequences

We next note that:

**Meta-lemma 15**: If the alphabet A of NT consists of only the numerals "0" and "1", then every converging NT-routine defines a Cauchy sequence[125] of rational numbers.

*Proof*: If $B$ is a converging NT-routine, then there is a positive integer $n_0$ such that, given any positive integer $n > n_0$, there is a minimum configuration number $n_{min}$ of output length $n$ such that:

(*i*)  the instantaneous tape description whose configuration number is $n_{min}+1$ is of the form $P_{n+1}q_jS_iQ$, where $P_{n+1}$ is a binary string of length $(n+1)$;

---

(*ii*)   the instantaneous tape description whose configuration number is $k > n_{min}+1$ is of the form $P_{n+1}P_mq_jS_iQ$ where $P_m$ is either the null string, or a binary string of finite, non-zero, length.

We thus have the Cauchy sequence $\{P'_0, P'_1, ..., P'_n, ...\}$, where $P'_n$ is obtained by prefacing a decimal point to the start of the string $P_n$, and is thus the binary form of a rational number between 0 and 1.¶

**Corollary 15.1**: If the alphabet A of an NT machine consists of only the numerals "0" and "1", a Cauchy sequence such as $\{P'_0, P'_1, ..., P'_n, ...\}$ may not be constructible if the instantaneous tape outputs of an NT machine define an oscillating, non-terminating, NT-routine.[126]

**Meta-lemma 16**: Every number-theoretic function that is effectively computable individually in the standard interpretation M of P defines a Cauchy sequence.

*Proof*: If a function $F(x_1, ..., x_n)$ of M is effectively computable individually in M, then it is effectively representable in P. Hence there is a P-formula $[R(y, x_1, ..., x_n)]$, with the free variables $[(y, x_1, ..., x_n)]$, such that, for any sequence of numerals $<[b], [a_1], ..., [a_n]>$ of P, $[(E!y)R(y, a_1, ..., a_n)]$ is P-provable.

The denumerable sequences $<[b], [a_1], ..., [a_n]>$ of P can be put into a 1-1 correspondence with the natural numbers. We thus have:

$k <-> <S_k>$

---

[126] In other words, even when a classical, non-halting, Turing machine corresponds to Turing's "circle-free" machine, its output does not necessarily define a real number, contrary to Turing's assertion in his seminal 1936 paper [Tu36]. Thus, Turing's "uncomputable" numbers may simply be the outputs of oscillating NT-routines. For instance, consider the circle-free machine that computes the non-terminating sequence of rationals $s_1, s_2, ...$, where $s_n$ = the first $n$ digits of $\sin(n)$.

for any given natural number $k$, where $<S_k>$ is a sequence of numerals $<[s_1], ..., [s_{n+1}]>$ of P. If we define the number-theoretic function $r_k$ such that $r_k = 0$ if, and only if, $[R(s_1, s_2, ..., s_{n+1})]$ is P-provable, and $r_k = 1$ otherwise, we obtain the Cauchy sequence:

$(0. \; r_1), (0. \; r_1 \; r_2), (0. \; r_1 \; r_2 \; r_3), ..., (0. \; r_1 \; r_2 \; r_3 ... \; r_k), ...$¶

**Corollary 16.1**: By Meta-lemma 11, it follows that every function that is effectively computable uniformly in M defines a Cauchy sequence.

**Definition 7**(*i*): We define the sequence:

$(0. \; r_1), (0. \; r_1 \; r_2), (0. \; r_1 \; r_2 \; r_3), ..., (0. \; r_1 \; r_2 \; r_3 ... \; r_k), ...$

as the *characteristic Cauchy sequence* of the function $F(x_1, ..., x_n)$ under the 1-1 correspondence $k <-> <S_k>$.

## 7.6 The P versus NP problem[127]

We now consider an argument that, by Corollary 1.2 to Meta-lemma 1, the definition[128] of the class P of polynomial-time languages in the P versus NP problem may not define a formal mathematical object.

**Definition 7**(*ii*): Let A be a finite alphabet (that is, a finite non-empty set) with at least two elements, and let A* be the set of finite strings over A. The language L(T),

---

[127] We follow the official description of the "The P versus NP Problem" provided by Professor Stephen Cook, University of Toronto, for the Clay Mathematical Institute. The downloadable pdf file is: <http://www.claymath.org/Millennium_Prize_Problems/P_vs_NP/_objects/Official_Problem_Description.pdf>

[128] This definition is based on the above description of the "The P versus NP Problem".

accepted[129] by the classical Turing machine whose description number is T[130], is defined by

L(T)= {w is in A* | T accepts w}

We denote by $t_T$(w) the number of steps in the computation of T on input w. For any given natural number $n$, we define the worst case run time of T as:

$T_T(n) = \max\{t_T(w) \mid w$ is in $A^n\}$

where $A^n$ is the set of all strings over A of length $n$. We say that T runs in polynomial time if there exists a $k$ such that, for all $n$, $T_T(n) =< n^k + k$. We then define the class P of languages by

P = {L | L = L(T) for some Turing machine T which runs in polynomial time}.

Is P a formal mathematical object?

Assuming that $T_T(n)$ is well-defined[131], we can define a number-theoretic relation $F$(T, $k$, $n$) that holds if, and only if, $T_T(n) =< n^k + k$. Clearly, if the alphabet A is of length $l$, there are a maximum of $l^n$ possible strings of length $n$. Hence, for a given T, $T_T(n)$ is effectively computable individually in the standard interpretation M of P.

It follows that, for a given T and $k$, and any given natural number $n$, there is also an individually terminating routine to determine whether the number-theoretic relation $F$(T, $k$, $n$) holds. However, there may not be any uniformly terminating routine to determine whether, for any given T and $k$, $F(\mathrm{T}, k, n)$ holds for any given natural number $n$.

---

[129] We note that T accepts w if, and only if, the computation of w terminates finitely in a special state designated as the accepting state of T.

[130] Following Turing [Tu36], we assume that each classical Turing machine can be assigned a unique machine description (natural) number.

[131] We note that, for $T_T(n)$ to be well-defined, we need to assume that, for any given $n$, and any given w in $A^n$, there is an effective method to determine whether the computation of T on input w is non-terminating. In some cases, there may be no such method.

In other words, in the absence of a specific proof, we cannot conclude that $F(\mathrm{T}, k, n)$ is effectively definable in the standard interpretation of standard PA; nor, in view of Meta-lemma 10, that $F(\mathrm{T}, k, n)$ necessarily defines a formal mathematical object.

Prima facie, there are also no grounds for concluding that the language L(T) defines a formal mathematical object; ipso facto, we cannot conclude, without specific proof, that the class, P, of languages defined by:

P = {L | (ET)(L = L(T) & (E$k$)(A$n$)$F$(T, $k$, $n$))},

is a formal mathematical object.

## 8.  Turing's computable and uncomputable numbers

In his seminal 1936 paper [Tu36], Turing did not explicitly define computable (*real*) numbers, but chose to variously describe their properties in terms of Turing machines as below. However, as we note in *Corollary 17.1*, Turing's assumption that every "circle-free" machine necessarily defines a Dedekind[132] real number is invalid.

(*a*)  *The "computable" numbers may be described briefly as the real numbers whose expressions as a decimal are calculable by finite means. ... A number is computable if its decimal can be written down by a machine.*

(*b*)  *If an a-machine prints two kinds of symbols, of which the first kind (called figures) consists entirely of 0 and 1 (the others being called symbols of the second kind), then the machine will be called a computing machine.*

---

[132] We follow Rudin's definition of a Dedekind real number ([Ru53], p9, def. 1.31).

*If the machine is supplied with a blank tape and set in motion, starting from the correct initial m-configuration, the subsequence of the symbols printed by it which are of the first kind will be called the* sequence computed by the machine.

*The real number whose expression as a binary decimal is obtained by prefacing this sequence by a decimal point is called the* number computed by the machine.

(*c*)   *If a computing machine never writes down more than a finite number of symbols of the first kind it will be called* circular. *Otherwise it is said to be* circle-free. *... A machine will be circular if it reaches a configuration from which there is no possible move, or if it goes on moving, and possibly printing symbols of the second kind, but cannot print any more symbols of the first kind.*

(*d*)   *A sequence is said to be computable if it can be computed by a circle-free machine. A number is computable if it differs by an integer from the number computed by a circle-free machine.*

(*e*)    *A computable sequence O is determined by a description (number) of a machine which computes O.*

(*f*)    *To each computable sequence there corresponds at least one description number, while to no description number does there correspond more than one computable sequence. The computable sequences and numbers are therefore enumerable. A number which is a description number of a circle-free machine will be called a satisfactory number.*

(*g*)    *The expression "there is a general process for determining ..." has been used throughout this section as equivalent to "there is a machine which will determine ...". This usage can be justified if and only if we can justify our definition of "computable". For each of these "general process" problems can be expressed as* *a problem concerning a general process for determining whether a given integer $n$ has a property $G(n)$ [e.g. $G(n)$ might mean "$n$ is satisfactory" or " $n$ is the Gödel representation of a provable formula"], and this is equivalent to computing a number whose $n$-th figure is 1 if $G(n)$ is true and 0 if it is false.*

(h)  *The computable numbers do not include all (in the ordinary sense) definable numbers. Let P be a sequence whose $n$-th figure is 1 or 0 according as $n$ is or is not satisfactory. It is an immediate consequence of the theorem of §8 that P is not computable.*

**Meta-lemma 17**: Turing's sequence $P$ in the above definition does not necessarily determine a Cauchy sequence.

*Proof*: By Turing's argument in §8 of his 1936 paper [Tu36], $P$ is not effectively computable uniformly in the standard interpretation M of any formal system of Arithmetic. Hence, either $P$ is effectively computable individually in M such that some sub-sequence $P_n$ of $P$ is self-terminating or, by the Corollary 15.1 to *Meta-lemma 15*, the sequence $P$ may define an oscillating NT-routine.¶

**Corollary 17.1**: Turing's argument in §8 of his 1936 paper [Tu36] does not establish the existence of a (classically uncomputable) Dedekind real number[133].

(j)  *It is (so far as we know at present) possible that any assigned number of figures of P can be calculated, but not by a uniform process. When sufficiently many figures of P have been calculated, an essentially new method is necessary in order to obtain more figures.*

---

[133] We note that classical theory does not insist on Dedekind's definition of a real number, but accepts Cantor's and Turing's non-constructive definition of a real number as any natural number that is followed by a period, and a non-terminating sequence of the integers 0, 1, ..., 9; further, as Hodges (op cit) remarks, it is classically accepted that "Turing showed it possible to give unambiguous definitions of real numbers which are not computable".

(In other words, *P* may be effectively computable individually, but not effectively computable uniformly, in any standard interpretation of a formal system of Arithmetic.)

## 9. Cantor's diagonal argument

A similar argument applies to Cantor's diagonal construction of an uncountable set of real numbers. Cantor's classical argument only establishes that there is no uniformly terminating routine for determining a 1-1 correspondence between the natural numbers and the real numbers between 0 and 1. However, there still may be an individually terminating routine for determining a 1-1 correspondence between the natural numbers and the real numbers between 0 and 1. In this case, any uniformly terminating routine that computes the sequence whose $n$-th digit is 1 if the $n$-th digit in the binary expression of the real number corresponding to $n$ is not 1, and 0 otherwise, may include a self-terminating computation of an NT machine, or it may correspond to an oscillating NT-routine. In such a case, the routine would not define a Dedekind real number. We thus have:

> **Meta-lemma 18**: We cannot conclude that the Dedekind real numbers are uncountable by Cantor's diagonal argument.

## 10. Constructivity and classical Quantum Mechanics

The introduction of constructive definitions of classical mathematical concepts may permit formal systems of standard Peano Arithmetic to model some of the more paradoxical concepts of Quantum Mechanics. For instance, consider the following argument:

> (*a*) Gödel has proved in his 1931 paper [Go31a] that there is a formula [$R(x)$] such that, for any given natural number $k$, [$R(k)$] is provable in any formal system of Arithmetic such as standard PA.

(*b*) Hence, for any given $k$, there is always some effective method for evaluating the arithmetic expression $R(k)$ in the standard interpretation M of PA.

(*c*) Gödel has also proved in the above paper that $[(Ax)R(x)]$ is not PA-provable.

(*d*) *Thesis*: There is no uniform effective method (algorithm/Turing machine) that can evaluate the arithmetic expression $R(n)$ for any given $n$ in M. (We note that (*d*) is a consequence of Corollary 9.2)

(*e*) Thus, $R(n)$ is individually computable, but not uniformly computable.

(*f*) *Theorem (provable by induction)*: For any given $k$, we can always define, using Gödel's Beta-functions, some effective method (algorithm/Turing machine) $T(k)$ that can compute $R(n)$ for all $n<k$, i.e. $T(k)$ terminates for all $n<k$, but it "loops" on input $k$. This follows since, clearly, all methods that evaluate $R(n)$ for all $n<k$ cannot be non-terminating on input $k$; this would imply that $R(k)$ is undefined, which would contradict (*b*).

(*g*) *Quantum interpretation*: The process of finding $T(k+1)$ can be corresponded, firstly, to the act of finding a suitable method of measuring the value $R(k)$ precisely, and, secondly, to the collapse of the wave function at $k$ as a result of the measurement; we then have the new "state" $T(k')$, which can evaluate the value of $R(n)$ for all $n<k'$, where $k<k'$, but not beyond!

(*h*) If, now, we have some law that determines the state $T(k')$ from the state $T(k)$ and the interaction at $k$, we have a deterministic interaction that is, nevertheless, absolutely unpredictable, where we may then define free will as absolute unpredictability. (We note that, if $k'>k+1$, we have a language that admits interactions that can leave the state $T(k')$ unchanged.)

Now we note that a counter-thesis to (*d*) would be:

(*i*) *Counter-Thesis*: There is some uniform effective method (algorithm/Turing machine) that can evaluate the arithmetic expression *R*(*x*) for any given *x*.

It follows from Gödel's reasoning in (1) that both (*d*) and (*i*) are effectively unverifiable, since they cannot be proved formally. We thus have two standard models of Peano Arithmetic - classical and constructive - that are mutually inconsistent. If we assume that both are consistent, the above argument indicates the interpretation that implies (*d*) may be the more suitable language for expressing concepts of classical Quantum Theory.

## References

[Go31a]  Gödel, Kurt. 1931. *On formally undecidable propositions of Principia Mathematica and related systems I*. In M. Davis (ed.). 1965. The Undecidable. Raven Press, New York.

[Go31b]  Gödel, Kurt. 1931. *On formally undecidable propositions of Principia Mathematica and related systems I*.
<*Web page*: http://home.ddc.net/ygg/etext/godel/index.htm>

[Ha47]  Hardy, G.H. 1947, 9[th] ed. Pure Mathematics. Cambridge, New York.

[La51]  Landau, E.G.H. 1951. Foundations of Analysis. Chelsea Publishing Co., New York.

[Me64]  Mendelson, Elliott. 1964. Introduction to Mathematical Logic. Van Norstrand, Princeton.

[Pe90]  Penrose, R. (1990, Vintage edition). The Emperor's New Mind: Concerning Computers, Minds and the Laws of Physics. Oxford University Press.

[Pe94]   Penrose, R. (1994). Shadows of the Mind: A Search for the Missing Science of Consciousness. Oxford University Press.

[Po01]   Podnieks, Karlis. 2001. Around Goedel's Theorem.

*<e-textbook*: http://www.ltn.lv/~podnieks/gt.html>

[Ru53]   Rudin, Walter. 1953. Principles of Mathematical Analysis. McGraw Hill, New York.

[Ti61]   Titchmarsh, E. C. 1961. The Theory of Functions. Oxford University Press.

[Tu36]   Turing, Alan. 1936. *On computable numbers, with an application to the Entscheidungsproblem.*

*<Web page*: http://www.abelard.org/turpap2/tp2-ie.asp - index>

(*Acknowledgements: My thanks to Professor Karlis Podnieks for suggesting, initially, the need for a comparison between Gödel's approach to his Theorem XI and that of the arguments underlying Meta-theorem 2. My thanks also to correspondents in various newsgroups who offered comments and counter-arguments that eventually led to the arguments of Meta-theorem 1, and to its consequences. My grateful thanks to Dr. Damjan Bojadziev for comments on the formatting and presentation. Author's e-mail: anandb@vsnl.com.*)

(*Updated: Saturday 10th May 2003 12:01:44 AM IST by re@alixcomsi.com*)

(*Created: Monday 5th Oct 2002 2:18:00 PM IST by re@alixcomsi.com*)